\documentstyle[12pt,twoside, epsf]{article}
\newcommand{\ZZ}{\mbox{$Z\!\!\! Z\!$}}	

 \newtheorem{th}{Theorem}
 \newtheorem{lem}{Lemma}
 \newtheorem{cor}{Corollary}
 \newtheorem{prop}{Proposition}

 \newtheorem{defn}{Definition}
 \newtheorem{conj}{Conjecture}

 \newtheorem{rem}{Remark}

\def\picill#1by#2(#3)
{\vbox to #2
{\hrule width #1 height 0pt depth 0pt
\vfill\epsffile{#3}}}

\let \ttorg \tt \def \tt{\ttorg \obeyspaces}

\begin{document}

\pagestyle{myheadings}

\date{}

\markboth{{\sc Kauffman \& Lambropoulou}}{{\sc Virtual Braids and the $L$--Move}}

\title{\bf Virtual Braids and the $L$--Move}

\author{\sc Louis H. Kauffman and Sofia Lambropoulou}

 \maketitle
  
 \thispagestyle{empty}

\section{Introduction}

In this paper we prove a Markov Theorem for the virtual braid group and for some analogs of this structure.
The virtual braid group is the natural companion to the category of virtual knots, just as the Artin braid 
group is  to classical knots and links. In classical knot theory the braid group gives a
fundamental algebraic structure associated  with knots. The Alexander Theorem tells us that every knot or
link can be isotoped to braid form. The capstone of this relationship is the  Markov Theorem, giving
necessary and sufficient conditions for two braids to close to the same link (where sameness of two links 
means that they are ambient isotopic). 
\bigbreak

The Markov Theorem in classical knot theory is not easy to prove. The Theorem was originally stated by
A.A.~Markov with three moves and then N. Weinberg reduced them to the known two  moves \cite{Ma,W}. The first 
complete proof is due to  J. Birman
\cite{Bi}. Other published  proofs are due to   D.~Bennequin \cite{Ben}, H. Morton \cite{Morton}, P. Traczyk
\cite{Traczyk} and S. Lambropoulou \cite{La,LR}. In this paper we shall follow the ``$L$--Move"  methods of
Lambropoulou. In  the $L$--move approach to the  Markov theorem, one gives a very simple uniform move that can
be applied anywhere in a braid to produce a braid with the same closure. This move, the
$L$--move, consists in cutting a strand of the braid and taking the top of the cut to the bottom of the braid
(entirely above or  entirely below the braid) and taking the bottom of the cut to the top of the braid
(uniformly above or below in correspondence with the choice for the other end of the cut). See Figure 15 for an 
 illustration of a classical $L$--move. One then proves that two braids have the same closure if  and only
if they are related by a sequence of $L$--moves. Once this
$L$--Move Markov Theorem is established, one can reformulate the result in various ways, including the more algebraic
classical Markov Theorem that uses conjugation and stabilization moves to relate braids with equivalent
closures.
\bigbreak

Up to now \cite{LR,HL,LR2} the $L$--moves were only used for proving analogues of the Markov theorem for
classical knots and links in $3$--manifolds (with or without boundary). 
Our approach to a Markov Theorem for virtual knots and links follows a similar strategy to the classical 
 case, but necessarily must  take into account properties of virtual knots and links that diverge from
the classical case. In particlular, we use $L$--moves that are purely virtual, as well as considering the effect
of allowed and forbidden moves of the virtual braids. The strategy for our project is to first give a specific
algorithm for converting a virtual link diagram to a virtual braid.  This algorithm is designed to be
compatible with the  $L$--moves. We prove that if two virtual diagrams are related to each other by a sequence
of virtual isotopy moves, then the corresponding braids are related by virtual $L$--moves and real conjugation. The exact
description of the $L$--moves for virtual braids is found in Section 2.2 (Definitions 2, 3, 4 and Figures 11, 12,
13). 
\smallbreak

The $L$--Move Markov Theorem for virtual braids is proved in Section 3 (see Theorem~2). Once the $L$--Move Theorem is
proved, it is a natural task to reformulate it in algebraic terms. In Section 4 we formulate and prove a local
algebraic Markov theorem for virtual braids (Theorem~3). This and the $L$--Move Markov Theorem for virtuals are the key
results of our paper. In Section 5 we recover  the Markov Theorem for virtual braids proved by  S. Kamada in
\cite{Ka}.  Such theorems are important for understanding the structure and classification of virtual knots and links. 
 The $L$--move approach provides a flexible conceptual center from which to
deduce many results. In particular, it would surely be quite difficult to compare our local algebraic
formulation of the Markov Theorem with that of Kamada without the fundamental $L$--move context. Our local algebraic
version of the Markov Theorem promises to be useful for formulating new invariants of virtual knots and links. 
\bigbreak

We conclude the paper with descriptions of variations of our Markov Theorem for other categories of braids,  
such as  flat virtual braids, welded braids and virtual unrestricted braids, in Section 6 (see Theorems 4, 5, 6, 7, 8).
For the case of welded braids our results coincide with the results of Kamada \cite{Ka}. Finally, in  Section 7
we describe the general pattern for obtaining quantum invariants via solutions to the Yang-Baxter equation and  Hecke
algebra type invariants of virtual links via braids.  These topics will be the subject of our future research.

\section{Virtual Knot Theory}

Virtual knot theory is an extension of classical diagrammatic knot theory. In this extension one 
adds  a {\em virtual crossing} (see Figure 1) that is neither an over-crossing
nor an under-crossing.  A virtual crossing is represented by two crossing arcs with a small circle
placed around the crossing point.  
\smallbreak

Virtual diagrams can be regarded as representatives for 
oriented Gauss codes (Gauss diagrams) \cite{VKT,GPV}. Some Gauss codes have planar
realizations, and these correspond to classical knot diagrams. Some codes do not have planar
realizations. An attempt to embed such a code in the plane leads to the production of virtual
crossings. 
\smallbreak

Another useful topological interpretation for virtual knot theory is in terms of embeddings of  links in
thickened surfaces, taken up to addition and subtraction of empty handles. Regard each  virtual crossing as a
shorthand for a detour of one of the arcs in the crossing through a 1-handle that has been attached to the
2-sphere of the original diagram (see Figure 1).   By interpreting each virtual crossing in this way, we
obtain an embedding of a collection of circles into a thickened surface  of genus the  number of virtual
crossings in the original diagram. See \cite{CKS,KamadaNS,VKT,DVK}. We say that two such surface embeddings
are {\em stably equivalent} if one can be obtained from another by isotopy in the thickened
surfaces,  homeomorphisms of the surfaces and the addition or subtraction of empty handles.  Then
we have the following Theorem \cite{DVK,CKS}: {\em Two virtual link diagrams are isotopic if
and only if their  correspondent  surface embeddings are stably equivalent.}  

\smallbreak

 A third way to make a topological interpretation of virtual knots and links is to form a
ribbon--neighborhood surface (sometimes called an {\it abstract link diagram} \cite{KamadaNS})  for a given
virtual knot or link, as  illustrated in Figure 1. In this Figure we show how a virtual trefoil knot (two
classical and one virtual crossing) has the classical  crossings represented as diagrammatic crossings in
disks, which are connected by ribbons, while the virtual crossings are represented by ribbons that pass over
one another without interacting. The abstract link diagram is shown embedded in three dimensional space, but
it is to be regarded without any particular embedding of the surface. Thus it can be represented with the
ribbons for the virtual crossings switched.  These abstract link diagrams give
the least surface embedding (with boundary) that can  represent a given link diagram.

\smallbreak

$$\vbox{\picill1.8inby2.3in(virtual1) }$$ 

\begin{center}
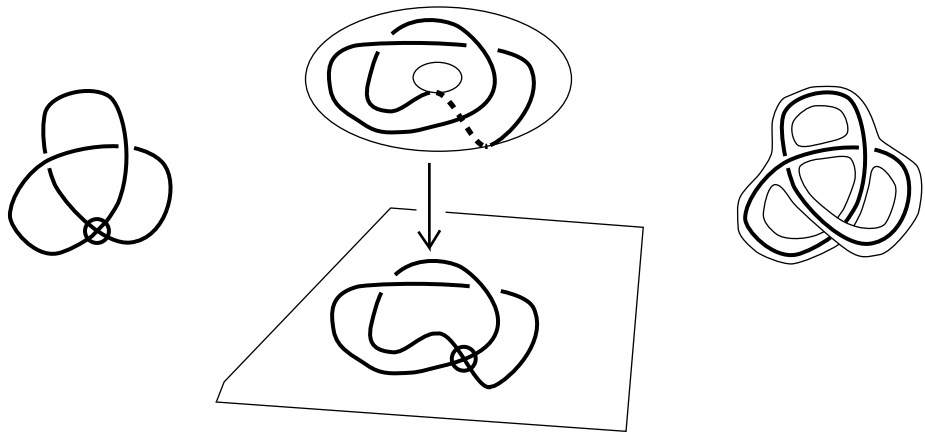
{\bf Figure 1 -- A  Virtual Trefoil and its Surface Realizations} 
\end{center}

Isotopy moves on virtual diagrams generalize the ordinary Reidemeister moves for classical knot and
link diagrams.  See Figure 2, where all variants of the moves should be considered. In this work,
virtual diagrams are always oriented, so the isotopy moves will be considered with all possible
choices of orientations. One can summarize the moves on virtual diagrams as follows:   The real
crossings interact with one another according to the {\it classical Reidemeister moves} (Part A of
Figure 2). Virtual crossings interact with one another by {\it virtual Reidemeister moves} (Part B of Figure 2). The
key move between virtual and classical crossings is shown in  Part C of Figure 2.  Here a consecutive sequence of two
virtual crossings can be moved across a single classical crossing. We will call it a
{\it special detour move}, because it is a special case  of the more general {\em detour move\/} indicated in Figure 3.
All these moves together with the {\it planar isotopy moves} (top left of Figure 2) generate an
equivalence relation in the set of virtual knot and link diagrams, called  {\it virtual
equivalence\/} or {\it virtual isotopy\/}.           

\bigbreak

$$ \picill5inby3.6in(virtual2)  $$

\begin{center} 
{ \bf Figure 2 --  Reidemeister Moves for Virtuals} 
\end{center}

 In the detour move, an arc in the 
diagram that contains a consecutive sequence of virtual crossings can be excised, and the arc
re-drawn, transversal to the rest of the  diagram (or itself), adding virtual crossings whenever
 intersections occur. See  Figure 3.  In fact, each of the moves in Parts B and C of Figure~2
can be regarded as special cases of  the detour move. By similar arguments as in
the classical Reidemeister Theorem, it follows that any detour move can be achieved by a finite
sequence of local steps, each one being a  Reidemeister move from Part B or C. 
A succinct description of virtual isotopy is that it is generated by classical Reidemeister moves
{\em and\/} the detour move.  

\bigbreak
 
$$ \picill3inby1.2in(virtual3)  $$

\begin{center} 
{ \bf Figure 3 -- The Detour Move} 
\end{center}

We note that a move analogous to a special detour move but with two 
real crossings and one virtual crossing is a {\it forbidden move\/} in virtual
knot theory.   There are two types of forbidden moves: One with an over arc,
denoted $F_1,$ and another with an under arc, denoted $F_2.$
 See \cite{VKT} for explanations and interpretations. Variants of the forbidden 
moves are illustrated in Figure 4.

\bigbreak
 
$$ \picill3inby.9in(virtual4)  $$

\begin{center} 
{ \bf Figure 4 -- The Forbidden Moves} 
\end{center}  

We know \cite{VKT,GPV} that classical knot theory embeds faithfully in virtual knot theory. That is, if two 
classical knots are equivalent through moves using virtual crossings, then they are equivalent as
classical knots via standard  Reidemeister moves.  With this approach, one can generalize many 
structures in classical knot theory to the virtual domain, and use the virtual knots to test the
limits of classical problems, such as  the question whether the Jones polynomial detects knots.  
Counter--examples to this conjecture exist in the virtual domain.  It is an
open problem whether some of these counter--examples are equivalent to  classical knots and links.

\section{Virtual Braids}

Just as classical knots and links can be represented by the closures of braids, so can  virtual
knots and links be represented by the closures of virtual braids \cite{SVKT, Ka, KL}.  A {\it virtual braid}
on $n$ strands is a braid on $n$ strands in the classical sense, which may also contain virtual
crossings. The closure of a virtual braid is formed by joining by simple
arcs the corresponding endpoints of the braid on its plane. Like virtual diagrams, a virtual braid can be
embedded in a ribbon surface. See Figure 5 for an example.

\bigbreak
 
$$ \picill2inby1.4in(virtual5)  $$

\begin{center} 
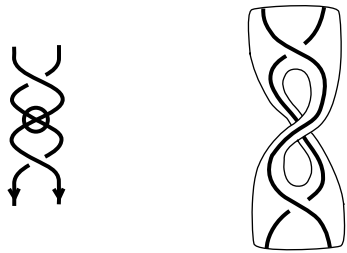
{ \bf Figure 5 -- A Virtual Braid and its Ribbon Surface Realization} 
\end{center}

The set of isotopy classes of
virtual braids on $n$ strands forms a group, the virtual braid group, denoted $VB_n,$ that can be
described by generators and relations, generalizing the generators and relations of the classical
braid group  \cite{SVKT}. This structure of virtual braids is worth study for its own sake. The virtual
braid group is an extension of the classical braid group by the symmetric group. See
\cite{VKT}, \cite{Bar}, \cite{KL}.  It is worth
remarking that classical  braids embed in virtual braids just as classical  links embed in virtual links.
This fact may be most easily  deduced from \cite{KUP}. 

\smallbreak

Virtual braids representing isotopic virtual links are related via a Markov--type virtual analogue.
In \cite{Ka}  S. Kamada proves a Markov Theorem for virtual braids, giving a set of moves on virtual
braids that generate the same equivalence classes as the virtual link types of their closures.
 For reference to previous work on virtual links and braids the reader should consult 
\cite{Bar,CKS,DK,GPV,HRK,KADOKAMI,KamadaNS,Ka,Ka1,VKT,SVKT,DVK,KL,KR,KUP,M1, M2,Satoh,TURAEV,V1}.

\subsection{Braiding  Virtual Diagrams}

It is easily seen that the classical Alexander Theorem \cite{A, Br} generalizes to
virtuals. 

\begin{th}{ \ Every 
(oriented) virtual link  can be represented by a virtual braid, whose closure is
isotopic to the original link.
}
\end{th}
Indeed, it is quite easy to braid a virtual diagram. In  \cite{KL} we gave, for example, 
a new braiding algorithm, which is applicable, in fact, to all the categories in which
braids are constructed. 
  The idea of that algorithm is very similar to the  braiding algorithm of Kamada \cite{Ka}, and it is the
following: we consider a virtual link diagram  arranged   in  general position  with respect to the
height function. We then rotate all crossings of the diagram on the plane, so that all arcs in the
crossings are oriented downwards. We leave the   down--arcs in place and eliminate the up--arcs,
producing instead braid strands. The elimination of an up--arc is described in Figure 9.

\bigbreak 

For the purposes of this paper, where we need to analyze  how the isotopy moves on diagrams
affect the final braids, we follow a different braiding process. 
\smallbreak 

\noindent {\bf Preparation for braiding.} Firstly, for simplicity and without loss
of generality, virtual link diagrams are assumed piecewise linear.  Working in the  piecewise linear
category gives rise to another `move': the {\it subdivision} of an arc into two smaller arcs, by marking
it with a point. The vertices and the local maxima and minima are subdividing points of a diagram. Subdivision of an
arc with no crossings can  be regarded as a degenerate case of the planar isotopy move.

\smallbreak 

Furthermore, virtual link diagrams lie on the plane, which is equipped with the top-to-bottom direction.
This  makes our set--up liable to certain conventions. For example, an oriented  virtual diagram contains
only {\it up--arcs} and {\it down--arcs} (no horizontal arcs). It contains no horizontally aligned
crossings, so as to have the crossings in the corresponding braid lying on different horizontal levels. 
Vertically aligned crossings or subdividing points are also not permitted, so as to avoid triple points when creating
new strands or pairs of braid strands with the same endpoints.  The above discussion  gives rise to the following
definition.

\begin{defn}{\rm \ 
 A virtual link diagram is said to be in  {\it general position} if it does not contain any
horizontal arcs and  no two  subdividing points or crossings are vertically or horizontally aligned, nor is a crossing
coincident with a maximum or a minimum.  }
\end{defn}

 Clearly, any virtual diagram can assume  general position by very small planar shifts. Note that,
the arcs or points or crossings that violate Definition~1 may not be close in the diagram. For example,
two aligned subdividing points may lie far away. The point is that the correcting shifts can
be applied on only one of them, so, in this sense these shifts can be assumed local.
 Moreover, when bringing a virtual diagram to general position we meet certain choices. For
example, in a parallel occurence of a maximum and a minimum, either one can occur first in the vertical
order.  Different choices amount to local shifts of crossings and subdividing points  with respect to
the horizontal or  the vertical direction. These local shifts shall be called {\it
direction sensitive moves}. 

The most interesting instances of such moves are  the {\it swing moves}. See Figure
6. A swing move avoids the coincidence of a maximum or minimum and a crossing, real or virtual. It turns
out  that adding the swing moves to our list of virtual  isotopy moves makes redundant certain instances of
Reidemeister moves involving horizontal arcs. For example, it is easily verified that an RII move with two horizontal
arcs can be produced by an RII move with two vertical arcs, two swing moves and changes of relative positions of
vertices. 

\smallbreak

$$\vbox{\picill3inby1in(virtual6)  }$$ 

\begin{center}
{ \bf Figure 6 -- The swing moves} 
\end{center}

It follows now easily that any two virtual diagrams in general position that correspond to isotopic
virtual diagrams will differ by the above `direction sensitive moves' and the  Reidemeister moves for
virtuals, all in general position. From now on, all diagrams will be assumed in general position.  

\smallbreak 
\noindent {\bf The braiding.} We are now ready to describe our braiding algorithm. The down--arcs will
stay in place while the  
  up--arcs shall be eliminated. Now, an up--arc will either be an arc of a crossing or it will be a
 free up--arc. We place each crossing containing one or two up--arcs in a small rectangular
box with diagonals the arcs of the crossing, the {\it crossing box}. A crossing box is assumed 
sufficiently narrow, so that the vertical zone it defines does not intersect the
zone of another crossing. The {\it free up--arcs} are  arcs joining the crossing boxes. We  first braid
the crossings containing an up--arc, one by one, according to the crossing charts of Figure 7. Except
for the local crossings shown in the illustrations, all other  crossings of the new braid strands with the rest
of the diagram are virtual. This is indicated abstractly by placing virtual crossings at the ends of the new
strands.  The result is a virtual tangle diagram. 

\smallbreak

$$\vbox{\picill8inby6.1in(virtual7)  }$$

\begin{center}
{ \bf Figure 7 -- The Braiding Chart for  Crossings} 
\end{center}

It is easy to verify that closing  the  corresponding braid strands of a braided crossing results in a virtual
tangle diagram isotopic to the starting one. In Figure 8 we illustrate this isotopy for one of the less
obvious cases.

\smallbreak

$$\vbox{\picill1.7inby3.2in(virtual8)  }$$

\begin{center}
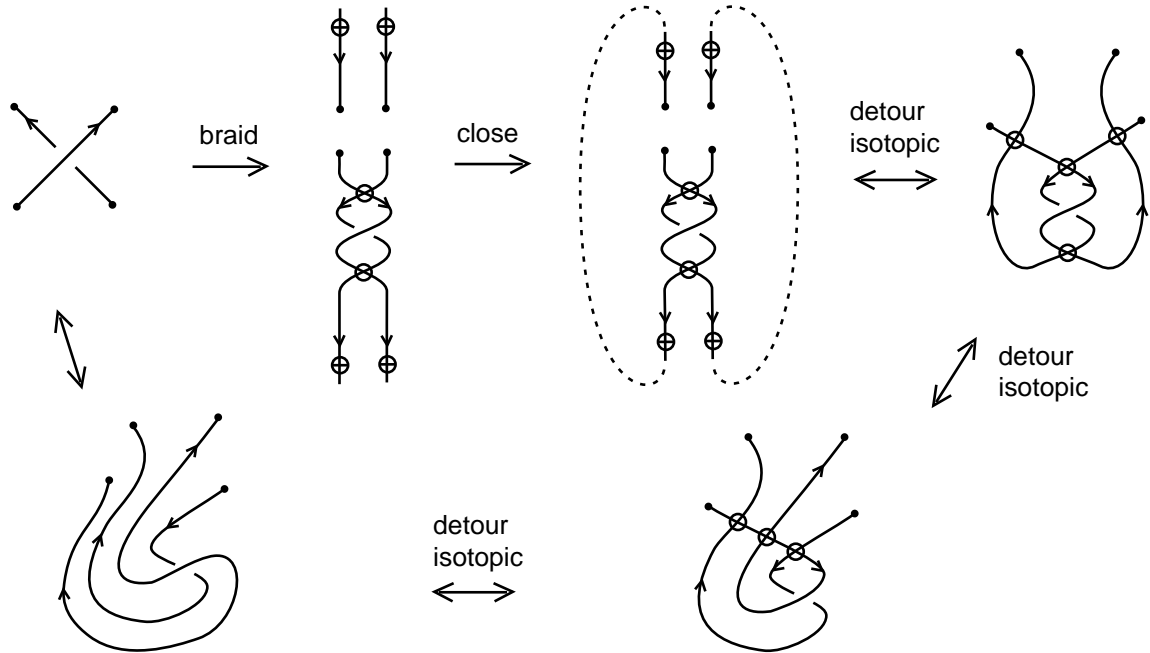
{ \bf Figure 8 -- The Closure of the Braiding of a Crossing} 
\end{center}

It remains to braid the free up--arcs. We braid a free up--arc by sliding it first across the 
right-angled triangle with hypotenuse the up--arc and with the right angle lying below it, so 
that {\it it  crosses virtually}  any other arcs of the original diagram that intersect the sliding triangle.
A grey curved arc is illustrated to this effect in Figure~9. We then  cut the vertical segment at a point and   we
pull the two ends, the upper upward and the lower downward, keeping them aligned, so that the two new braid
strands cross any other part of the  diagram  {\it only virtually}. This is indicated in the illustrations by
the virtual crossings on the final braid strands. We also care that the horizontal  arc slopes slightly
downwards, so that there is no conflict with not permitting horizontal arcs in Definition~1. Note  that 
the prior elimination of crossings may cause vertical strands to cross virtually the free up--arc.  This is
not an obstacle for braiding it, since --by the detour move--  the arc can  slide virtually across these
strands (see grey strands in Figure 9). 
\smallbreak

 In the end we created  a pair of corresponding braid strands and we have one  up--arc less. 
 Note that  joining  the two  corresponding braid strands yields a virtual tangle diagram obviously
isotopic to the starting one, since from the free up--arc we created a stretched loop around the
braid axis, which is detour isotopic to the  arc. The braiding of a free up--arc is a {\it basic
braiding move}.

\smallbreak 

$$\vbox{\picill4inby3in(virtual9)  }$$ 

\begin{center}
{ \bf Figure 9 --  The  Basic Braiding Move}
\end{center}

 After completing all braidings we obtain an open virtual 
braid, the closure of which is  an oriented virtual link diagram  isotopic to the original one.
 The braiding algorithm given above will braid any virtual diagram and, thus, it proves  Theorem~1. 
\hfill$\Box$  

\smallbreak

$$\vbox{\picill1.3inby2.3in(virtual10)  }$$

\begin{center}
{ \bf Figure 10 --  An Example of Braiding} 
\end{center}

\begin{rem} \rm Because of the narrow zone condition for the crossings (see the beginning of the  braiding
discussion) the braidings of the crossings are independent, so their order is irrelevant. Moreover, because of
the braid detour move, it does not make any difference in which order we braid the free up--arcs. In fact,
we could even braid any number of them before completing the braidings of the crossings.
\end{rem}

\begin{rem} \rm  The braidings of the crossings are also based on the basic braiding move. Using
this, it is easy to verify that, if in the instances of the braiding chart we replace   each arc by a
number of parallel arcs with the same orientation and the same crossings, the resulting braids are
$L$--equivalent to the ones we would obtain if we braid one by one the single crossings in the
formation, according to the chart. This remark can save us from creating unnecessary extra braid
strands.
\end{rem}

 The set--up of our virtual braiding resembles the one in \cite{LR} for classical
links, but only to the extent that we consider piecewise linear diagrams on the plane, which is equipped with
the top-to-bottom direction, and that the basic braiding move looks similar. With the forbidden moves in the theory, the
choices needed here are  completely different from the ones made in the classical set--up. For example, we are
forced  to braid a crossing of two up--arcs as one entity, not its arcs one by one.  (Braiding
crossings as rigid entities can, obviously, be applied also in the classical set--up for braiding knots and
links.)  In the classical set--up, braiding an over up--arc corresponds to pulling the new pair of braid strands over the
rest of the diagram. Here it has to be always virtually. Another technical difference is that, in the classical set--up
it was important to ensure that the sliding triangles have no intersections with other parts of the diagram. Here  this
assumption is not needed.

\subsection{The $L$-equivalence for Virtual Braids}

As in classical knot theory, the next consideration after the braiding is to characterize virtual
braids that induce, via closure, isotopic virtual links. In this section we describe an
equivalence relation between virtual braids, the {\it $L_v$--equivalence}. For this purpose we need to
recall and generalize to the virtual setting the {\it $L$--moves} between braids to {\it virtual $L$--moves}, abbreviated
to {\it $L_v$--moves}. The $L$--move (see Definition 5) was introduced in \cite{La,LR}, where it was used among other
things to  prove the `one--move Markov theorem' for classical oriented links (cf.Theorem 2.3 in \cite{LR}), replacing the
two well-known moves of the Markov equivalence: the {\it stabilization} that introduces a crossing at the bottom right of
a braid and {\it conjugation} that conjugates a braid by a crossing.

\begin{defn}{\rm \ 
A  {\it basic $L_v$--move}  on a virtual braid, 
consists in cutting an arc of the  braid open and pulling the upper cutpoint downward and
the lower  upward, so as to create a new pair of braid strands with corresponding endpoints
(on the vertical line of the cutpoint), and such that both strands  cross entirely  {\it virtually}
 with the rest of the braid. (In abstract illustrations this is indicated by placing virtual crossings on
the border of the braid box.) 
 }
\end{defn}

By a small braid isotopy that does not  change the relative positions of endpoints,
a basic $L_v$--move can be equivalently seen as introducing an in--box virtual crossing to a virtual braid,
which faces either the {\it right} or the {\it left} side of the braid. If we want to
emphasize the existence of the virtual crossing, we will say {\it virtual $L_v$--move}, abbreviated to {\it $vL_v$--move}.
In  Figure 11 we give abstract illustrations. See also Figure 16 for a concrete example.

\smallbreak

$$\vbox{\picill4inby1.7in(virtual11)  }$$ 

\begin{center}
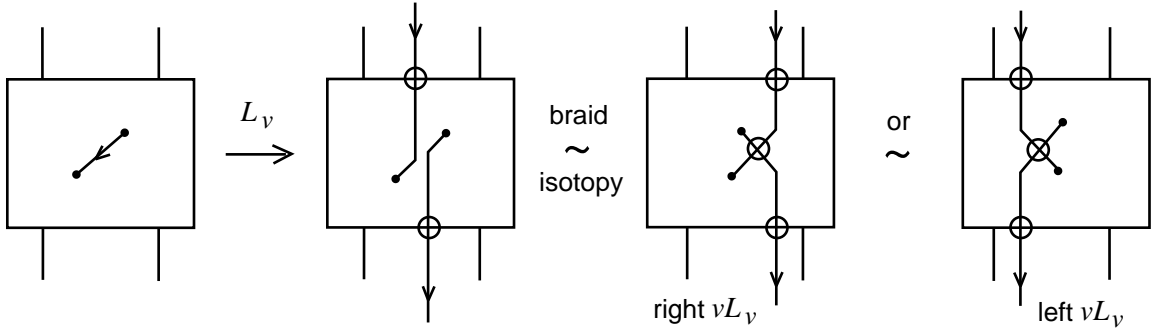
{ \bf Figure 11 -- A Basic $L_v$--move and the two $vL_v$--moves} 
\end{center}

 Note  that in the closure of a basic $L_v$--move or a $vL_v$--move the detoured loop contracts to a kink. This kink could
also be created by a real crossing, positive or negative. So we define:

\begin{defn}{\rm \ 
A  {\it real $L_v$--move}, abbreviated to  {\it  $+L_v$--move} or {\it  $-L_v$--move}, is a virtual 
$L$--move that introduces  a real in--box crossing (positive or negative) on a virtual braid, and it can face either the
{\it right} or the {\it left} side of the braid. See  Figure 12 for abstract illustrations.}
\end{defn}

\smallbreak

$$\vbox{\picill4.2inby1.6in(virtual12)  }$$ 

\begin{center}
{ \bf Figure 12 --  Left and Right  Real $L_v$--moves} 
\end{center}

If the crossing of the kink is virtual, then, in the presence of the forbidden moves, there
is another possibility for an $L_v$--move  on the braid level, which uses another arc of the braid, the `thread'. So we
have: 

\begin{defn}{\rm \ 
A  {\it threaded $L_v$--move} on a virtual braid is a virtual  $L$--move with a virtual
crossing in which, before pulling open the little up--arc of the kink, we perform a Reidemeister II
move with real crossings, using another arc of the braid, the {\it thread}. See Figure 13. There are two possibilities:
an  {\it over--threaded $L_v$--move} and an {\it under--threaded $L_v$--move}, depending on whether we pull the kink 
over or under the thread,  both with the variants {\it right} and {\it left}.   
}
\end{defn}

\smallbreak

$$\vbox{\picill4.2inby1.7in(virtual13)  }$$ 

\begin{center}
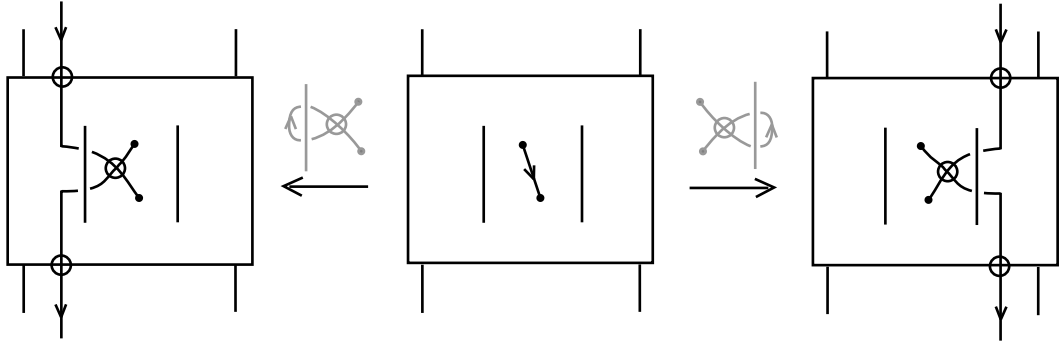
{ \bf Figure 13 -- Left and Right Under--Threaded $L_v$--moves} 
\end{center}

Note that a threaded $L_v$--move cannot be simplified in the braid. 
 If the crossing of the kink were real, then, using a braid RIII move  
with the thread, the move would reduce  to a real $L_v$--move. Similarly, if the
forbidden moves were allowed, a threaded $L_v$--move would reduce to a  $vL_v$--move.

\begin{rem}{\rm \  As with a braiding move, the effect of a  virtual $L$--move, basic, real or threaded, is to 
stretch (and cut open) an arc of the braid around the braid axis using the detour move, after twisting it and possibly
after threading it. Conversely, such a move between virtual braids gives rise to isotopic closures,
since the virtual $L$--moves shrink locally to kinks (grey diagrams in Figures 12 and 13). }
\end{rem}

 Conceivably, the `threading' of a virtual $L$--move could involve a
sequence of threads and Reidemeister II moves with over, under or virtual crossings, as Figure 14
suggests. The presence of the forbidden moves does not allow for simplifications on the braid level.  We show
later that such {\it multi--threaded $L_v$--moves}  follow from the threaded $L_v$--moves, up to real conjugation. 

\smallbreak

$$\vbox{\picill4inby1.65in(virtual14)  }$$ 

\begin{center}
{ \bf Figure 14 -- A Right Multi--Threaded $L_v$--move} 
\end{center}

We finally introduce the notion of a classical $L$--move, adapted to our set--up. 

\begin{defn}{\rm \ A {\it classical $L_{over}$--move} resp.  {\it $L_{under}$--move } on a virtual braid
consists in  cutting an arc of the virtual braid open and pulling the two ends, so as to create a new
pair of braid strands, which run both {\it entirely over}  resp. {\it entirely under}  the rest of the
braid, and such that the closures of the virtual braids before and after the move are isotopic. See
Figure~15 for abstract illustrations.  A classical  $L$--move may also introduce an in--box
crossing, which may be positive, negative or virtual,  or it may even involve a thread. }
\end{defn}

In order that a classical $L$--move between virtual braids is {\it allowed}, in the sense that it gives rise to
isotopic virtual links upon closure, it is required that the virtual braid  has no virtual crossings on the
entire vertical zone either to the left or to the right of the new strands of the $L$--move. We then perform the isotopy 
on the side with no virtual crossings. 
 We show later that the allowed  $L$--moves can be expressed in terms of 
$L_v$--moves and real conjugation. It was the classical $L$--moves that were introduced  in \cite{LR}, and they
replaced the two equivalence moves of the classical Markov theorem. Clearly, in the classical set--up these
moves are always allowed, while the presence of forbidden moves can preclude them in the virtual setting.

\smallbreak

$$\vbox{\picill4inby1.6in(virtual15)  }$$ 

\begin{center}
{ \bf Figure 15 -- The Allowed Classical  $L$--moves} 
\end{center}

 In Figure~16  we illustrate an example of various types of $L$--moves taking place at the same point
of a virtual  braid. 

\smallbreak

$$\vbox{\picill5inby1.5in(virtual16)  }$$ 

\begin{center}
{ \bf Figure 16 -- A Concrete Example of Introducing $L$--moves} 
\end{center}

\section {The $L$--move Markov Theorem for Virtual Braids}

It is clear that different choices when applying the braiding algorithm  as well as local 
isotopy changes on the diagram level may result in  different virtual braids. 
In this section we show that {\it real conjugation} (that is, conjugation by a real crossing)
and some variations of the $L_v$--moves (recall Definitions 2, 3, 4) capture and reflect on the braid
level all instances of isotopy between virtual links.

\begin{th}[$L$--move Markov Theorem for virtuals] Two oriented  virtual links are isotopic if and
only if any two corresponding virtual braids differ by virtual braid isotopy and a finite sequence of the
following moves or their inverses:
\begin{itemize}
\item[(i)]Real conjugation 
\item[(ii)]Right virtual $L_v$--moves  
\item[(iii)]Right real $L_v$--moves  
\item[(iv)]Right and left under--threaded $L_v$--moves.
\end{itemize} 
\end{th}

\begin{defn}{\rm \ Moves  (i), (ii), (iii), (iv) together  with virtual
braid isotopy generate an equivalence relation in the set of virtual braids, the {\it $L$--equivalence}, used in the
statement of Theorem~2. }
\end{defn}

Note that in the statement of Theorem~2 we do not use virtual conjugation, basic $L_v$--moves, left virtual or real
$L_v$--moves, allowed classical $L$--moves, over--threaded $L_v$--moves (right or left) and multi--threaded
$L_v$--moves. In the next lemmas we show that all these moves (except for the left real $L_v$--moves) follow from the
$L$--equivalence. We shall then use them freely in the proof of Theorem~2. The proof that left real $L_v$--moves follow
from the $L$--equivalence shall be given at the end of the proof of the Theorem (Lemma 9).

\begin{lem}{ \ Virtual conjugation can be realized by a sequence of basic and virtual $L_v$--moves.
 }
\end{lem}

\noindent {\it Proof.} \ The proof is an adaptation for virtual conjugation of a similar proof of 
R.~H\"{a}ring-Oldenburg for classical braids and real conjugation \cite{HL}. In Figure~17 we start with a
virtual braid conjugated by $v_i.$ After performing an appropriate basic $L_v$--move and braid isotopy, and 
after undoing  another virtual $L_v$--move we end up with the original braid. $\hfill \Box$  

\bigbreak

Note that the `trick' of Figure 17 would not work in the case of real conjugation. In fact, we conjecture  that
real conjugation cannot be generated by virtual $L$--moves.
 
\smallbreak

$$\vbox{\picill8inby1.6in(virtual17)  }$$

\begin{center}
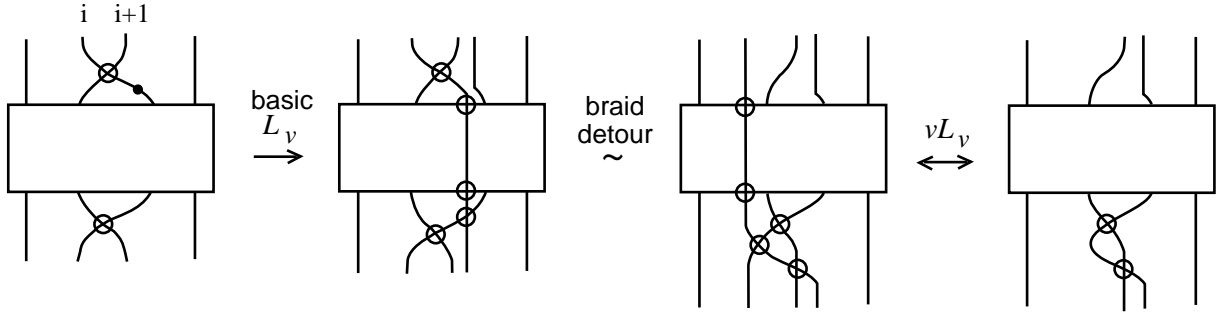
{ \bf Fig. 17 -- Conjugation by $v_i$ is a composition of $L_v$--moves}
\end{center}
\vspace{3mm}

\begin{lem}{ Basic and left virtual $L_v$--moves follow from right virtual $L_v$--moves and braid isotopy. } 
\end{lem}

\noindent {\it Proof.} \ The proof is illustrated in Figure 18.  $\hfill \Box$   

\smallbreak

$$\vbox{\picill4inby2in(virtual18)  }$$ 

\begin{center}
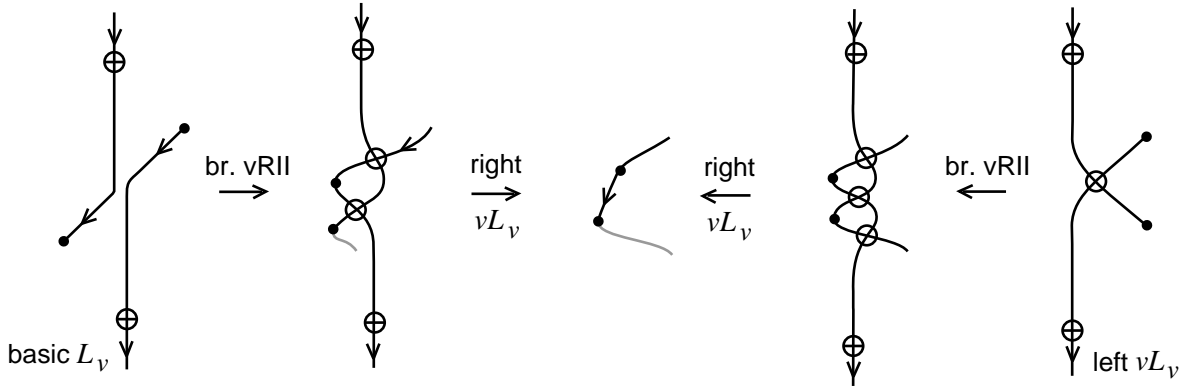
{ \bf Figure 18 --  Basic and Left Virtual $L_v$--moves as Right $vL_v$--moves} 
\end{center}

It is easy to see that an allowed classical $L$--move reduces, up to real conjugation and classical braid relations, to
a right or left real $L_v$--move at the extreme right or left of the braid box. See Figure~37 and the discussion after
Remark~7.

\bigbreak

We shall now prove a key lemma about `in--box exchange moves'.

\begin{defn}{\rm \ An {\it in--box exchange move} is a move between virtual braids as illustrated in Figure~19 between
the first two or the last two pictures, together with the two variants with facing the
left (obtained by reflecting the diagrams in a vertical axis).
}
\end{defn}

\smallbreak

$$\vbox{\picill5inby1.7in(virtual19)  }$$ 

\begin{center}
{ \bf Figure 19 -- In--box Exchange Moves} 
\end{center}

\begin{lem}{ \ The in--box exchange moves follow from  $L_v$--moves and real conjugation.
 }
\end{lem} 

\noindent {\it Proof.} \ The proof is illustrated in Figure 20. For the second step we point out that the real
conjugation we do here can be carried out just above and below the middle box, for the following reason: since the $i$th
and $(i+1)$st strands cross the top and bottom braid boxes virtually, any crossing in their vertical zone can be braid
detoured away. So, this vertical zone will only contain parts of other strands (drawn in grey), crossing the $i$th and
$(i+1)$st strand virtually. Then, the real crossing formed by these two strands can be braid detoured to the top, get
real conjugated to the bottom and then pass in the same manner to the region above the bottom box. $\hfill \Box$

\smallbreak

$$\vbox{\picill8inby2.1in(virtual20)  }$$ 

\begin{center}
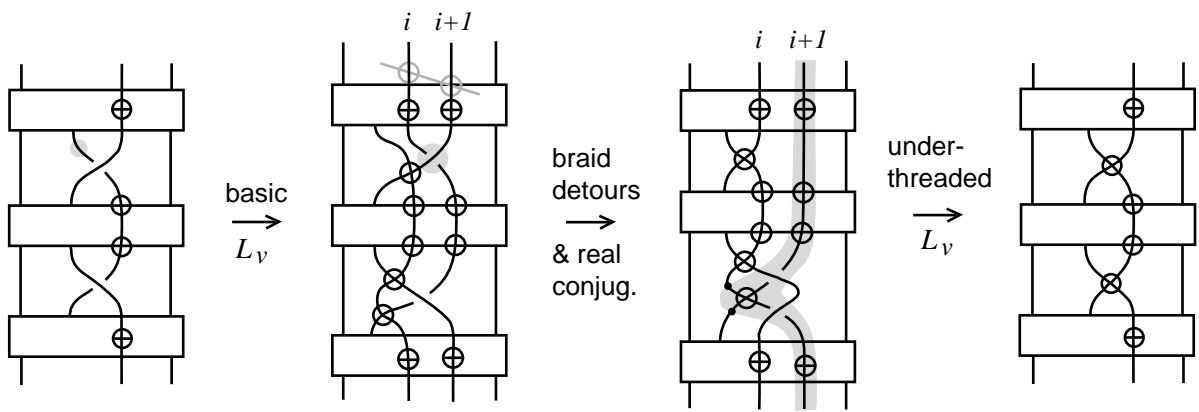
{ \bf Figure 20 -- An In--box Exchange Move via $L_v$--moves and Real Conjugation} 
\end{center}

The next lemma shows that in the $L$--equivalence we only need, indeed, one type of threaded $L_v$--moves, say the
under--threaded (left and right), recall Figure~13.

\begin{lem}{ \ The over--threaded $L_v$--moves follow from the $L$--equivalence moves of Definition~6.
 }
\end{lem} 

\noindent {\it Proof.} \ Lemma~3 is the key. Indeed, as illustrated in Figure~21, a right over--threaded
$L_v$--move gives rise to an in--box exchange move of the same type as the one in Figure~20. So, applying Lemma~3
involves only $L$--equivalence moves. Similarly, the in--box exchange move facing the left (with a top
negative real crossing) involves only a virtual $L_v$--move and a left under--threaded $L_v$--move. Thus, a left 
over--threaded $L_v$--move follows also from $L$--equivalence moves.  $\hfill \Box$

\smallbreak

$$\vbox{\picill8inby2in(virtual21)  }$$ 

\begin{center}
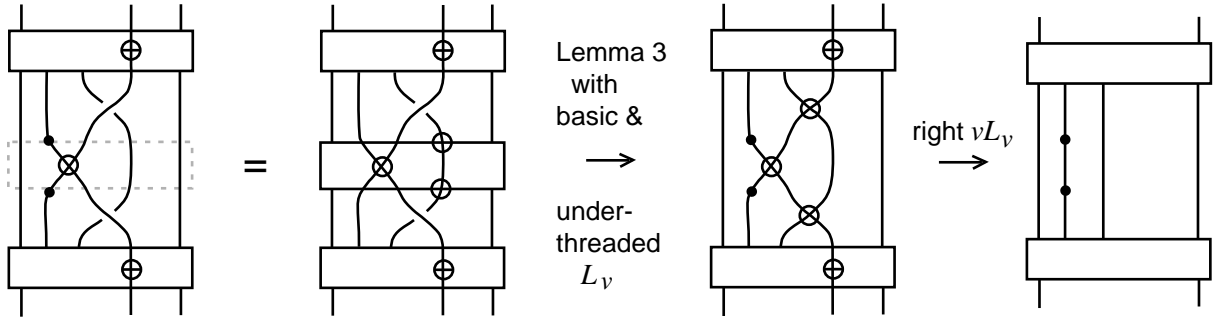
{ \bf Figure 21 -- An Overthreaded $L_v$--move via the $L$--equivalence Moves} 
\end{center}

As a result of Lemmas 3 and 4 and their proofs we have the following.

\begin{cor}{ \  The in--box exchange moves follow from the $L$--equivalence moves. }
\end{cor} 

\noindent {\it Proof.} \ Indeed, when in the proof of Lemma~3 we reach an over--threaded $L_v$--move, we apply the
process in Figure~20 with just one virtual crossing in the middle box. $\hfill \Box$

\begin{lem}{ \ The multi--threaded $L_v$--moves are consequences of the $L$--equivalence moves.
 }
\end{lem} 

\noindent {\it Proof.} \ Notice first that all
threads can be assumed real, as virtual threads can be braid--detoured away around the virtual crossing of
the move. In Figure 22 we illustrate the last stage of the proof. We assume any number of real threads inside the middle
braid box, instead of just one illustrated here. Then, by Lemma~3 we exchange the two real crossings of the
outer thread with two virtual ones and we braid--detour away the virtual thread. We proceed like this until we are
left with one real thread. If in the application of Lemma~3 or in the last step an over--threaded $L_v$--move is created,
apply Lemma~4 and Corollary~1.
 $\hfill \Box$

\smallbreak

$$\vbox{\picill8inby1.7in(virtual22)  }$$ 

\begin{center}
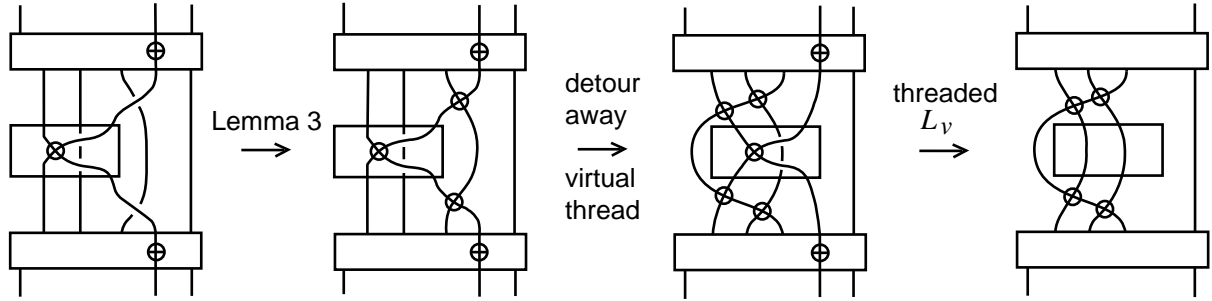
{ \bf Figure 22 -- A Multi--threaded $L_v$--move Follows from the $L$--equivalence Moves} 
\end{center}

\begin{rem}{\rm \ The in--box exchange moves of Definition 7 generalize the virtual exchange moves defined  by
S.~Kamada in \cite{Ka}, which he used in formulating and proving a Markov type theorem for virtual braids. See Section~5 
for details.
 }
\end{rem}

\noindent {\bf Proof of Theorem 2.} \   
Clearly, $L$--equivalent braids have isotopic closures. We have to show the converse. 

\smallbreak 
The proof  splits into two parts: the {\it technical part} and the {\it isotopy part}. In the
technical part we compare virtual braids resulting from  different choices made on a given virtual diagram 
during the braiding process. The isotopy part consists in comparing virtual braids corresponding to  virtual
diagrams that are  related either by different choices made when bringing a diagram to general position
(recall Definition~1) or by  the virtual isotopy moves. 
\smallbreak 

We first discuss  the technical part. Since our braiding is quite rigid, the only choices made during the
braiding process are the subdividing points and the order of the braiding moves.  The order of the braiding
moves is irrelevant, according to Remark~1.  Subdividing points  are  needed
for marking the up--arcs and the crossing boxes.  Assume now that our diagram is equipped with a choice of subdividing
points. In order to compare it to a different choice of subdividing points we need the following lemma.

\begin{lem}{If we add to an up--arc  an extra subdividing point, the corresponding braids differ by basic 
$L_v$--moves.} 
\end{lem}

\noindent {\it Proof.} \ Assume first that the up--arc is a free up--arc.  
Without loss of generality we have eliminated all other free up--arcs and crossings containing up--arcs
of the  diagram except for the up--arc in question and its subdivided replacement. We complete the
braiding by eliminating the up--arc. In Figure 23 we let $P$ be the new subdividing point of
the up--arc and  $P'$  its  projection on the horizontal arc (slightly sloping downwards) created
by the braiding. We perform a  basic  $L_v$--move at $P'$ and, by a small braid planar isotopy, we
obtain the braid that would result from the original diagram with the subdividing point $P$ included. 

\smallbreak

$$\vbox{\picill4inby1.5in(virtual23)  }$$ 

\begin{center}
{ \bf Figure 23 -- The Proof of Lemma 6} 
\end{center}

If, now, the  up--arc  is an arc inside a crossing box, then we create a similar smaller box inside
the original, using the new subdividing point, and we complete the braiding of the new formation.
Again, we will find that the corresponding braids differ by two or four basic $L_v$--moves,
depending on whether the crossing contains one or two up--arcs.
$\hfill \Box$

\begin{cor}{Given any two subdivisions $S_1$ and $S_2$ of a virtual diagram, the corresponding braids 
are $L$--equivalent.}
\end{cor}
 
\noindent Indeed, consider the subdivision $S_1 \bigcup S_2$, which is a common refinement of $S_1$ and
$S_2$, and apply repeatedly Lemma 1 to $S_1$ and to $S_2$.

\bigbreak 

We proceed now with the isotopy part of the proof of Theorem~2. The choices we have when bringing a
virtual diagram to general position are related to the  direction sensitive
moves  (recall discussion after Definition~1). 
 These,  as well as the virtual isotopy moves, are all local. Thus, given two virtual diagrams that differ
by such a move, we may assume that they  have both been   braided everywhere, except for the arcs and
crossings inside the  regions of the local move. After  completing the braiding, we compare two
virtual braids, which are identical except for the effect of the move on each.  In the figures that
follow we  focus only on the local moves and their braidings,  dropping the abstract  box.

\begin{lem}{Virtual diagrams in general position that differ by direction sensitive moves correspond to
virtual braids that differ by basic and virtual $L_v$--moves.} 
\end{lem}

\noindent {\it Proof.} \ Repairing a horizontal arc corresponds to a planar isotopy
move. If the arc in the move is an up--arc,  the move boils down to subdivision of an up--arc (Lemma~6),
basically because subdivision can be seen as a degenerate case of planar isotopy. (We refer the reader to 
\cite{LR} for  details.) In Figure  24 we check  planar isotopy  in the case of a down--arc.  

\smallbreak

$$\vbox{\picill4inby2in(virtual24)  }$$ 

\begin{center}
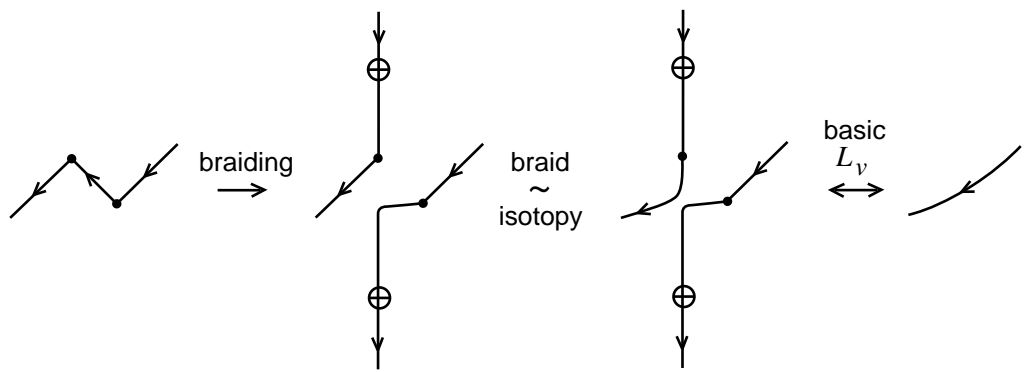
{ \bf Figure 24 -- Checking one Case of Planar Isotopy} 
\end{center}

Changes of relative heights of crossings or subdividing points yield --up to virtual braid relations--
the same virtual braids. Also, vertical alignment of crossings or subdividing points can be 
repaired by local sidewise shifts. In Figure 25 we illustrate a case of vertical alignment
 of two subdividing points and its braided resolutions. We only show the two up--arcs
containing the subdividing points (everything else is already braided), the alignment of which is
indicated by a dotted line. Note that, up to  a virtual braid RII move, the two braids are conjugates by
a virtual crossing.
 All other cases of vertical alignment are based on the same idea, possibly involving conjugation by more
than one virtual crossing.

\smallbreak

$$\vbox{\picill4inby2.4in(virtual25)  }$$ 

\begin{center}
{ \bf Figure 25 -- An Instance of Vertical Alignment} 
\end{center}

We shall now check the swing moves. There are various cases, depending on the orientation, the type of
crossing and the minimum/maximum. The ones with a virtual crossing are very easy to check. In Figures 26
and 27 we check two cases with a real crossing. 

\smallbreak

$$\vbox{\picill4inby2in(virtual26)  }$$ 

\begin{center}
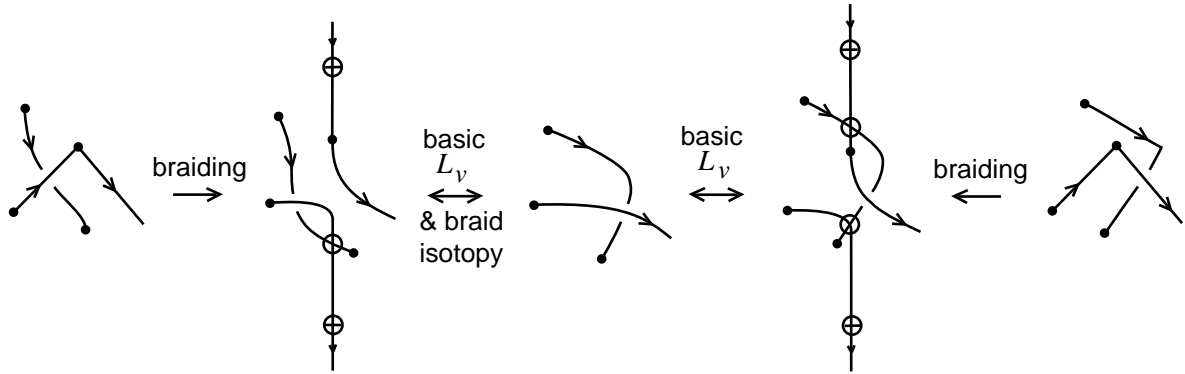
{ \bf Figure 26 -- The First Case of the Swing Moves} 
\end{center}

\smallbreak

$$\vbox{\picill4inby2.3in(virtual27a)  }$$ 

$$\vbox{\picill4inby2.2in(virtual27b)  }$$ 

\begin{center}
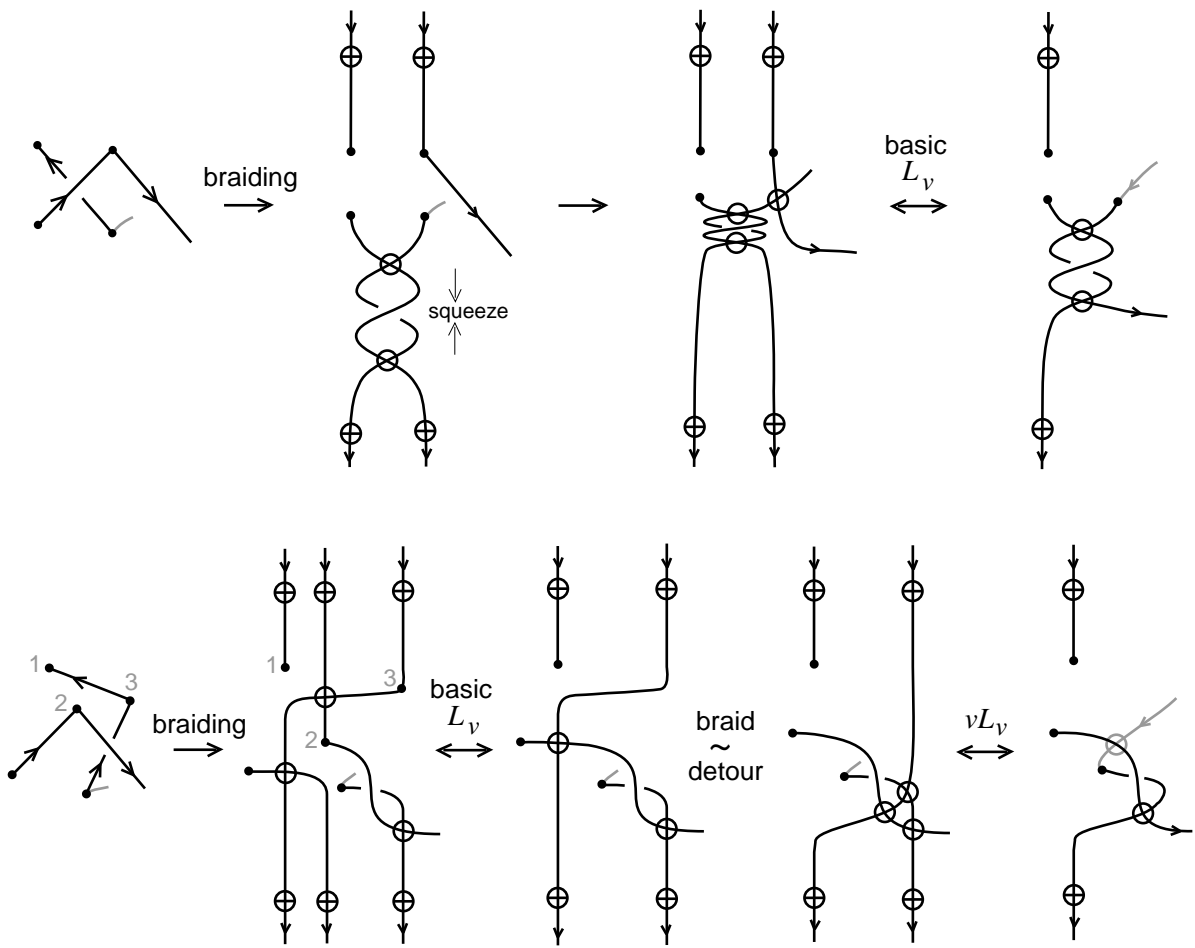
{ \bf Figure 27 -- The Second Case of the Swing Moves} 
\end{center}

In many parts of Figure 27 we have drawn in grey the continuation of an arc. In the last instance this is needed for 
 comparing the final braids of the two sides of the move.  The key point here is that  this grey arc
 is part of the braiding of an up--arc, so its crossing with our braid diagram will be virtual. Recall our
assumption, that in the regions of the local moves there are no other crossings of the original diagram.
Note, finally, that if this continuation arc was pointing to the left, so it would be in both sides of the
move and, again, the two final braids would agree. 
$\hfill \Box$  

\bigbreak 

We shall now check the virtual isotopy moves. Indeed we have:

\begin{lem}{Virtual diagrams that differ by virtual isotopy moves correspond to virtual braids that differ by braid
isotopy, $L_v$--moves and real conjugation.} 
\end{lem}

\noindent {\it Proof.} \ We discuss first the RII moves. An RII move with
two down--arcs follows immediately from virtual braid isotopy.  From all  cases of RII moves (with virtual/real
crossings, different orientations) the vertical ones checked in Figures 28 and 29 are the most interesting ones. For
example, an RII move placed horizontally follows from two swing moves and a vertical RII move.
 In Figure 28 we check a reverse real RII move. Note here that a threaded
$L_v$--move is involved. If this were an over--threaded $L_v$--move, we could apply Lemma~4 in order to use only
an under--threaded one.

\smallbreak

$$\vbox{\picill4inby2.3in(virtual28)  }$$ 

\begin{center}
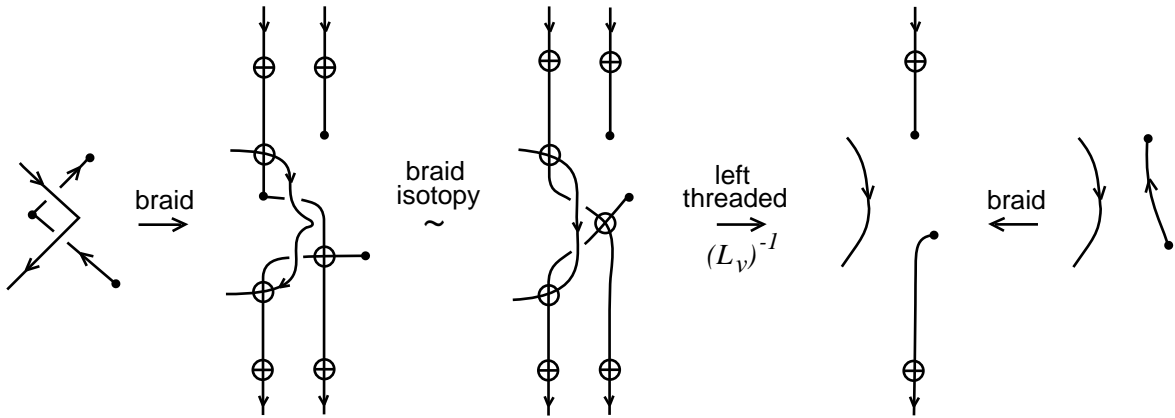
{ \bf Figure 28 -- The First Case of an RII Move} 
\end{center}

In Figure 29 we check a real RII move with two up--arcs. Here it is real conjugation that
will play the main role. Again, the braiding has been done for the rest of the diagram, and parts of
it are indicated in grey. The braiding algorithm ensures that  there are no other real crossings in
the final braid lying in the narrow vertical zone, which is created after the completion of the
braiding. So, apart from the real crossings indicated, all other crossings in this vertical zone will be virtual
(drawn in grey), created by the new braid strands. As in the proof of Lemma~3, 
this means that the old strands act as channels for the real crossings 
to reach the top and the bottom of the braid, hence to be available for conjugation.

\smallbreak 

$$\vbox{\picill8inby6in(virtual29)  }$$ 

\begin{center}
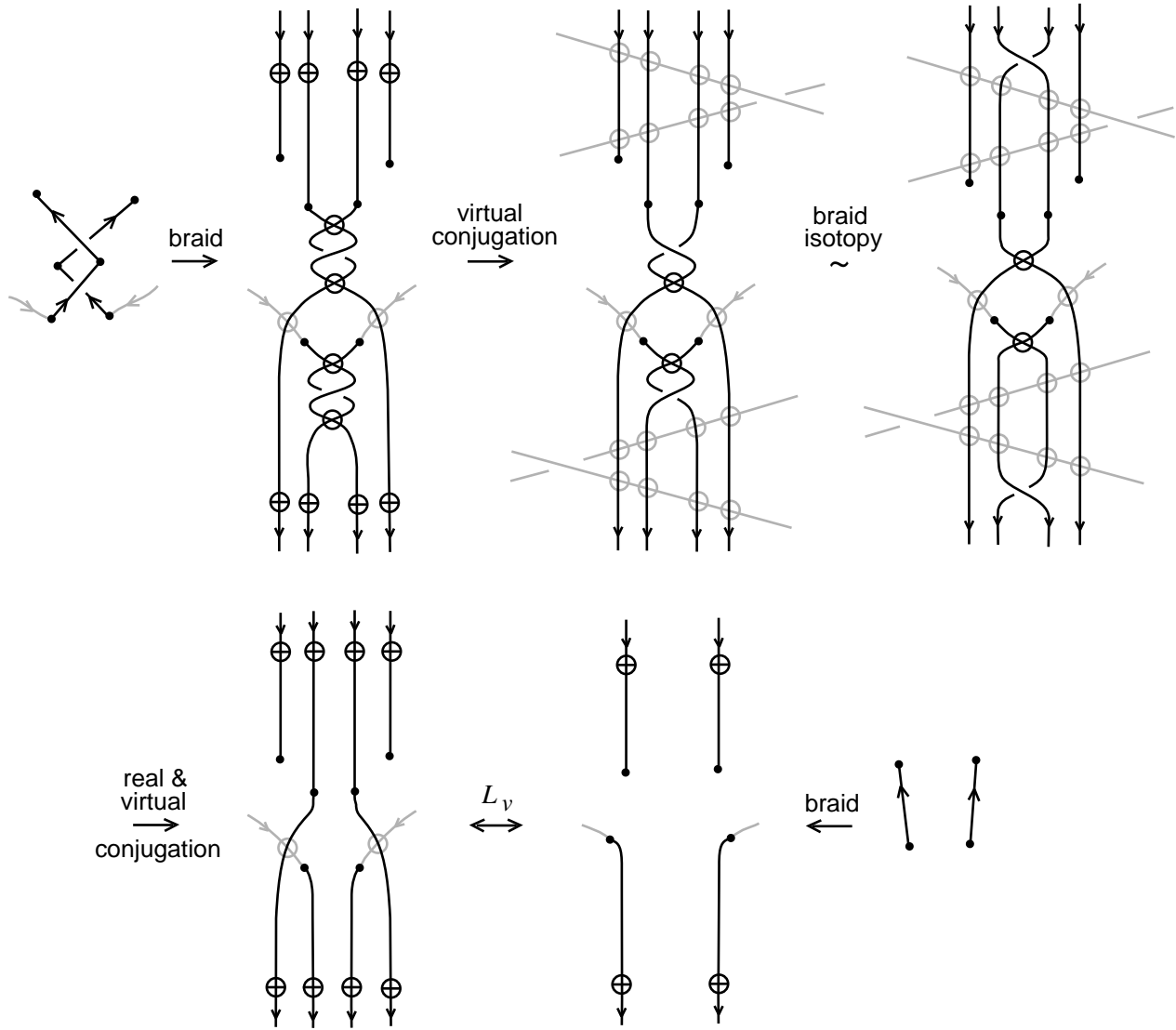
{ \bf Figure 29 -- The Second Case of an RII Move} 
\end{center}

We shall now check  the  RIII type moves. These include the classical RIII moves, the virtual RIII moves and the
special detour moves. Note that all RIII type moves  with three down--arcs are, in fact,  braid  relations. Consider now
an RIII move with one up--arc and two down--arcs. Using a well--known trick we can perform the move using RII moves
(which are already checked) and an RIII move with three down--arcs. See Figure~30. The same trick applies to virtual
and special detour moves, but in some cases of special detour moves 
 we may have to also use the swing moves. In Figure 31 we demonstrate the most interesting
case. If an RIII type move involves two up--arcs and one down--arc we apply the same trick to reduce to the case
of one up--arc. Similarly, an RIII type move with three up--arcs reduces to the previous cases. $\hfill \Box$  

\smallbreak

$$\vbox{\picill4inby1.1in(virtual30)  }$$ 

\begin{center}
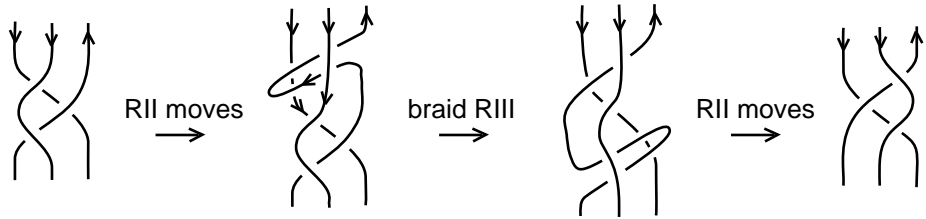
{ \bf Figure 30 -- A Real RIII Move with One Up--arc} 
\end{center}

\smallbreak

$$\vbox{\picill8inby1in(virtual31)  }$$ 

\begin{center}
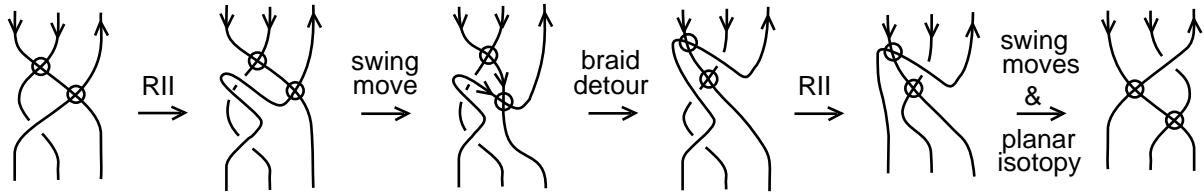
{ \bf Figure 31 -- A Special Detour Move with One Up--arc} 
\end{center}

We shall finally check  the  RI type moves.
Virtual RI moves on the diagram clearly give rise to $vL_v$--moves on the braid level. Braiding
an RI move with downward orientations on the crossing will clearly give rise to a right real $L_v$--move, if the kink
faces the right, and to a left real $L_v$--move, if the kink faces the left. Braiding
an RI move with upward orientations on the crossing will give rise to a right real $L_v$--move, if the kink faces the
left (see Figure 32), and to a left real $L_v$--move, if the kink faces the right. $\hfill \Box$  

\smallbreak

$$\vbox{\picill4inby2.2in(virtual32)  }$$ 

\begin{center}
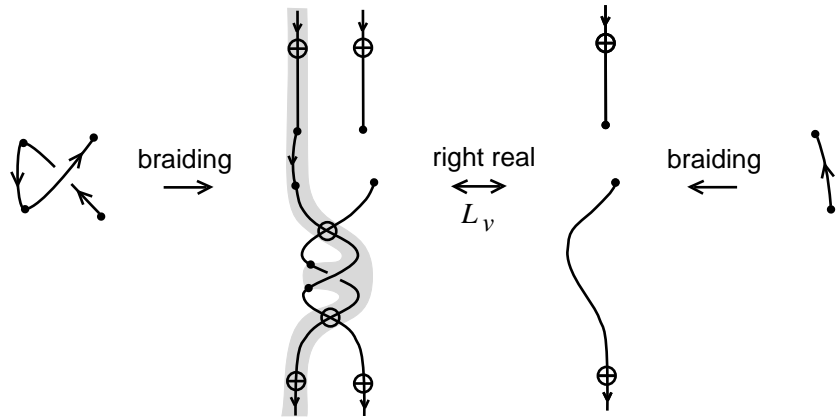
{ \bf Figure 32 -- An RI Move giving rise to a Right Real $L_v$--move} 
\end{center}

Lemma~9 below completes the proof that all $L_v$--move follow from $L$--equivalence moves.

\begin{lem}{ A left real $L_v$--move  can be performed by a sequence of $L$--equivalence moves. Consequently, an RI move
giving rise to a left real $L_v$--move corresponds to a sequence of $L$--equivalence moves.} 
\end{lem}

\noindent {\it Proof.} \  For the proof we employ the Whitney trick (compare
\cite{Ka}). In Figure 33 we start with a virtual diagram $K_1,$ which is almost the closure of a braid $B,$
except that it contains a kink with a real crossing introduced in $B$. So, $K_1$ opens to a virtual braid $B_1,$
which contains a left real $L_v$--move. On 
$K_1$ we introduce a second kink and we perform a sequence of isotopy moves that undo the kink we started
with. At the same time we register at each step the difference that every isotopy move makes on the braid
level. The final diagram $K_9$ is, then, the closure of the starting braid $B.$  So, we went from
$B_1,$ containing a left real $L_v$--move, to $B$  with the  $L_v$--move removed, via a sequence of
$L$--equivalent braids (Definition~6). $\hfill \Box$   

\smallbreak

$$\vbox{\picill8inby4.9in(virtual33)  }$$ 

\begin{center}
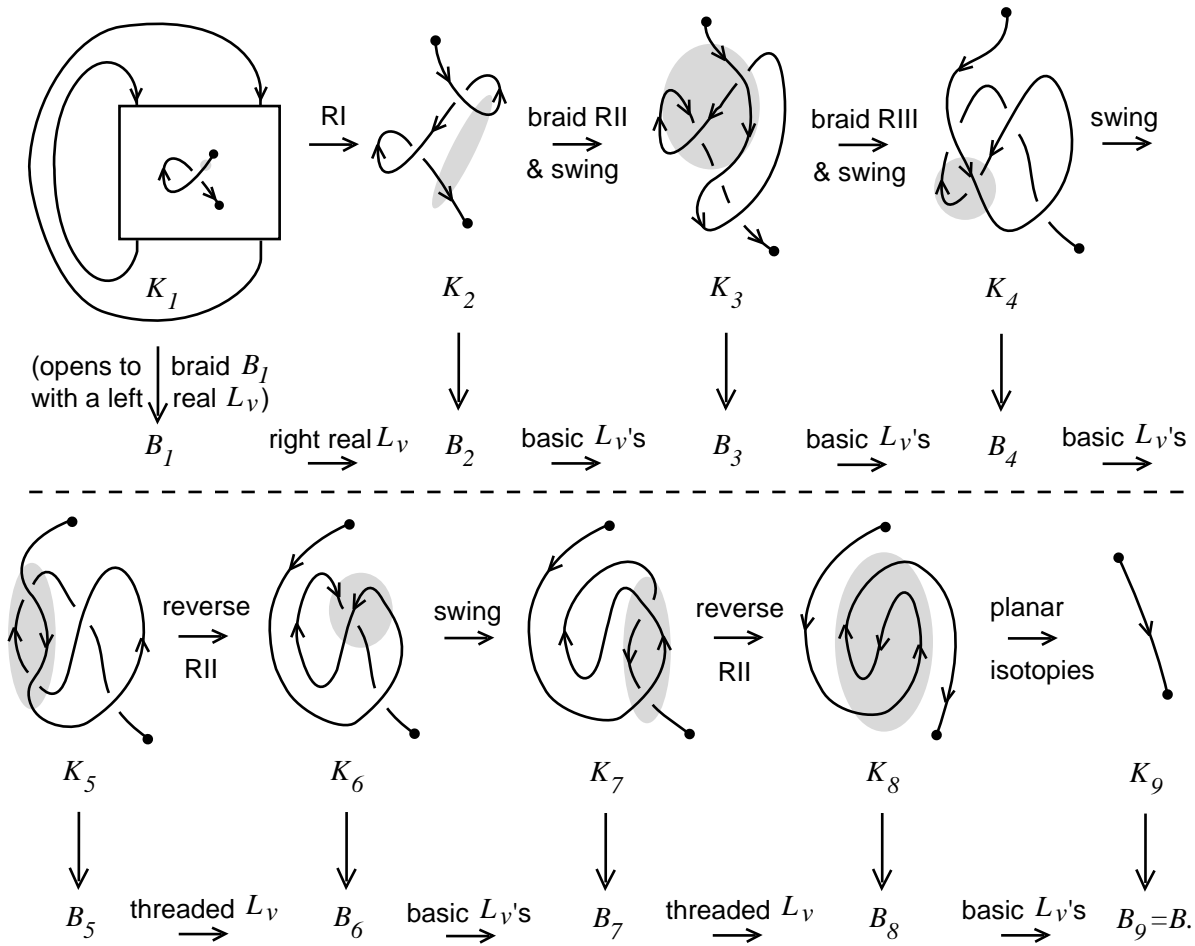
{ \bf Figure 33 -- The Whitney Trick for the Left Real $L_v$--moves} 
\end{center}

By Lemmas 1, 2, 4,  7, 8, 9 and  by Corollaries 2 and 3 the proof of Theorem 2 is now concluded. 
$\hfill \Box$  

\begin{rem} \rm As far as the proof of Theorem~2 is concerned, the reverse real RII isotopy moves are the only cases
where the threaded $L_v$--moves appear on the braid level. Moreover, real conjugation is needed in the proof of Lemma~4
and in a real RII with two up--arcs. Finally, the real $L_v$--moves appear only in the real RI cases.
\end{rem} 

\begin{conj} \rm {\it Real conjugation is not a consequence of the $L_v$--moves.} In other words, it should be
possible to construct a virtual braid invariant  that will not distinguish  $L_v$--move equivalent virtual
braids, but will distinguish virtual braids that differ by real conjugation. As the simplest possible
puzzle, try to show that there is no sequence of $L_v$--moves connecting the  pair of equivalent
braids shown in Figure~34.
\end{conj}

\smallbreak

$$\vbox{\picill1inby1.2in(virtual34)  }$$ 

\begin{center}
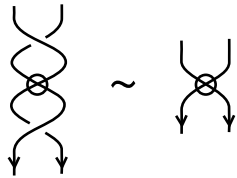
{ \bf Figure 34 -- The Simplest Pair of Real Conjugates} 
\end{center}

\vspace{.4in}

\section {Algebraic Markov Equivalence for Virtual Braids}

In this section  we reformulate and sharpen the statement of Theorem~2 by giving an equivalent list of local
algebraic moves in the virtual braid groups. More precisely, let  $VB_{n}$ denote the virtual braid group on
$n$ strands and let $\sigma_i, v_i$ be its generating classical and virtual crossings. The $\sigma_i$'s
satisfy the relations of the classical braid group and the $v_i$'s satisfy the relations of the permutation
group. The characteristic relation in $VB_{n}$  is the {\it special detour move} relating both:
 
$$ \begin{array}{cccl} 
v_{i} \sigma_{i+1} v_{i} & = & v_{i+1} \sigma_i v_{i+1}. &   \\ 
\end{array}$$

\noindent The group $VB_{n}$ embedds naturally into $VB_{n+1}$ by adding one identity strand at the right of
the braid. So, it makes sense to define $VB_{\infty} :=
\bigcup_{n=1}^{\infty} VB_{n}$,  the disjoint union of all virtual braid groups. We can now state our result.

\begin{th}[Algebraic Markov Theorem for virtuals] Two oriented  virtual links are isotopic if and
only if any two corresponding virtual braids differ by a finite sequence of braid relations in $VB_{\infty}$ and the
following moves or their inverses:

\begin{itemize}
\item[(i)]Virtual and real conjugation:    \ \ \ \ \ \ \ \ \ \ \ \  $ v_i \alpha v_i \sim \alpha \sim 
{\sigma_i}^{-1}\alpha \sigma_i $
\item[(ii)]Right virtual and real stabilization:  \ \ \ \ \  $\alpha v_n \sim \alpha
\sim \alpha \sigma_n^{\pm 1}$ 
\item[(iii)]Algebraic right under--threading:  \ \  $\alpha \sim \alpha \sigma_n^{-1} v_{n-1} \sigma_n^{+1} $
\item[(iv)]Algebraic left under--threading:  \ \ \ \ $\alpha \sim  \alpha v_n v_{n-1} \sigma_{n-1}^{+1} v_n
\sigma_{n-1}^{-1} v_{n-1} v_n $,
\end{itemize} 

\noindent where $\alpha,  v_i, \sigma_i \in VB_n$ and  $v_n, \sigma_n \in VB_{n+1}$ (see Figure 35).
\end{th}

\smallbreak

$$\vbox{\picill3inby1.2in(virtual35a)  }$$ 

$$\vbox{\picill5inby1.9in(virtual35b)  }$$ 

\begin{center}
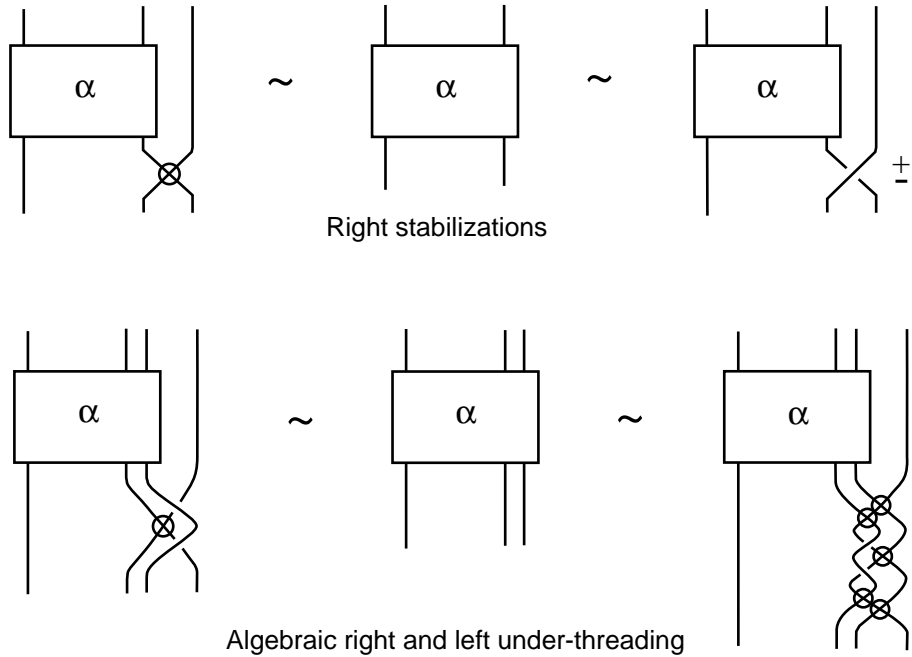
{ \bf Figure 35 -- The Moves (ii), (iii) and (iv) of Theorem 3} 
\end{center}

\begin{rem} \rm    Given $b$ in $VB_{n}$ let $i(b)$ denote the element of $VB_{n+1}$ obtained by adding one to 
the index of every generating element in $b$ (compare \cite{Ka}). In other words, $i(b)$ is obtained by adding a
single identity strand to the left of $b.$ We also regard $b$ as an element of $VB_{n+1}$  by adding a strand
on the right, but take this inclusion for granted, with no extra notation. In the above notation, a left
under--threaded $L_v$--move pulled to the bottom left side of the braid will have the  algebraic
expression:
$ \alpha \sim  i(\alpha) \sigma_1^{\pm 1}  v_2 \sigma_1^{\mp 1}$ (see Figure 36).

\smallbreak

$$ \vbox{\picill4inby1.4in(virtual36)  }$$

\begin{center} 
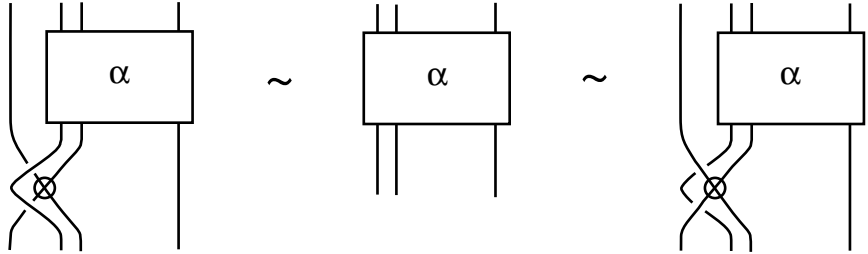
{ \bf Figure 36 -- Bottom left under/over threadings: $ \alpha \sim  i(\alpha) \sigma_1^{\pm 1}  v_2
\sigma_1^{\mp 1}$} 
\end{center}  

\end{rem} 

\noindent {\bf Proof of Theorem 3.} \ The algebraic moves of Theorem~3 follow immediately from the moves of Theorem~2 by
braid detouring to the right and by conjugation in $VB_{\infty}$. For example, in Figure~37 we illustrate how to bring a
right real $L_v$--move to the right end of the braid. In order to derive the algebraic left under--threaded moves: we
first bring a left under--threaded $L_v$--move to the bottom left of the braid by conjugation, and then we braid detour
to the right and apply virtual conjugation.  $\hfill \Box  $

\smallbreak

$$\vbox{\picill8inby1.5in(virtual37)  }$$ 

\begin{center}
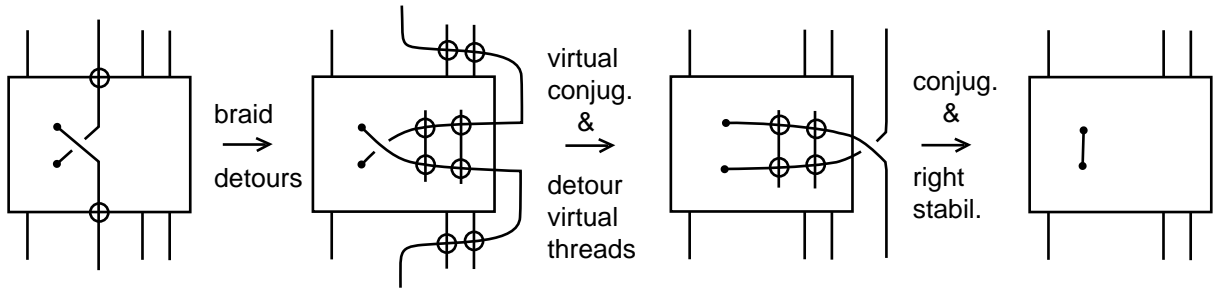
{ \bf Figure 37 -- Right Real $L_v$--move derived from  Right Stabilization} 
\end{center}

\begin{rem} \rm   By the braid conjugation, moves (ii), (iii), (iv) of Theorem~3 could be
equally given with the local algebraic part in between two braids. For example: 
$$\alpha \beta \sim \alpha \sigma_n^{-1} v_{n-1} \sigma_n^{+1} \beta.$$
\end{rem} 

Finally, we should point out that the proof in Figure 37 can be also adapted  to the case of  allowed
classical $L$--moves, namely pulling to the right or left, depending on which side is free of virtual
crossings. Here, the conjugation for pulling aside is real and agrees with the type (over/under) of the
classical $L$--move. Once out of the braid box, we have a real stabilization move.

\section {Kamada's Markov Theorem for Virtual Braids}

In this section we present  Kamada's Markov Theorem for virtual braids \cite{Ka} and we show that our
Theorem~3 is  equivalent to the Theorem of Kamada. With the inclusion of braids of Remark~6, S.~Kamada proved the
following:
 \bigbreak

\noindent {\bf Theorem} {\it (S.~Kamada \cite{Ka}) \ Two virtual braids $b$ and $b'$ have isotopic closures if and only
if they are related to one another through a finite sequence of braid relations in $VB_{\infty}$ and the following moves:

\begin{itemize}
\item[1.]{\em conjugation} if $b'$ is the conjugation of $b \in VB_n$ by an element of $VB_n$,
\item[2.]{\em right stabilization move} if $b'$  is $b\sigma_{n}$ or 
$b\sigma_{n}^{-1}$ or $bv_{n} \in VB_{n+1},$  for $b \in VB_n,$
\item[3.]{\em right exchange move} if they belong to one of the following patterns, for $b_1, b_2 \in VB_n$ 
$$\{ b_{1}\sigma_{n}^{-1}b_{2}\sigma_{n}, \, b_{1}v_{n}b_{2}v_{n} \,\},$$
\item[4.]{\em left exchange move} if they belong to one of the following patterns, for $b_1, b_2 \in VB_n$
$$\{ i(b_{1})\sigma_{1}^{-1}i(b_{2})\sigma_{1}, \, i(b_{1})v_{1}i(b_{2})v_{1} \,\}.$$
\end{itemize} 
}

\bigbreak
\noindent In Figure 38 we illustrate the braids for the right and left exchange moves. It is clear from the
Figure that the corresponding braids for these moves have equivalent closures.

\smallbreak

$$ \vbox{\picill4inby1.6in(virtual38)  }$$

\begin{center} 
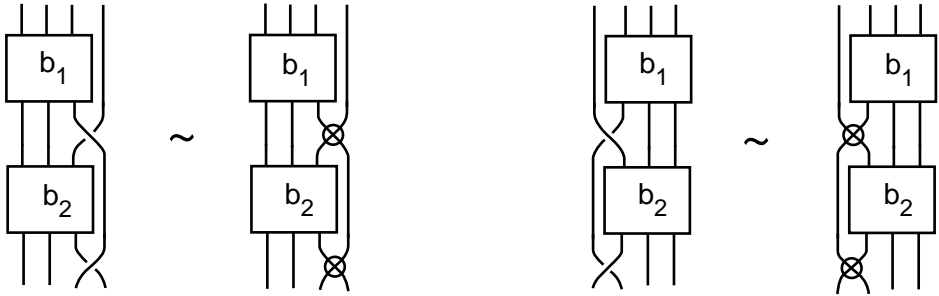
{ \bf Figure 38 -- The Right and Left Exchange Moves of Kamada} 
\end{center}

\begin{prop}  The moves of the Kamada Theorem follow from the moves of Theorem~3. Conversely, the moves
of Theorem~3 can be realized via the moves of  Kamada. 
\end{prop} 

\noindent {\em Proof.} \ The first two moves coincide with moves (i) and (ii) of Theorem~3. Further, Kamada's exchange
moves are special cases of in--box exchange moves (recall Definition~7), so, by Lemma~3 and by Theorems~2 and 3, they
follow from the moves of Theorem~3. 

Consider now an algebraic right under--threaded move. The one side of the move is a special case of one side of the
exchange move, where the second braid box contains only the virtual crossing $v_{n-1}.$ Perform the exchange
move to change the thread to a virtual one. Braid detour it away and apply the right virtual stabilization. This
brings us to the other side of the threaded move. Note now that Kamada's left exchange move is equivalent, up to
conjugation, to a similar left exchange move with the opposite crossings. Let us call that one an `under left exchange
move'. For realizing an algebraic left under--threaded move: conjugate it first to the bottom left of the braid (as in
Figure~36) and realize this move via an under left exchange move. Finally, conjugate the result back to the bottom right
of the braid. $\hfill \Box  $

\section {The Markov Theorem for Flat Virtuals and Welded Links }

In this section we give the analogues of Theorems 2 and  3 for flat virtuals, welded links
and virtual unrestricted links. Each category is interesting on its own right and has been studied by various 
authors. In \cite{KL} we gave reduced presentations for the corresponding braid groups.

\subsection {Flat Virtuals}

Every classical
knot or link diagram can be regarded as an immersion of cirlces in the plane with extra under/over structure 
at the double points.  If we take the diagram without this extra structure, it is the shadow of some link in
three dimensional space, but the weaving of that link is not specified. We call these shadow crossings {\it
flat crossings}. Clearly, if one is allowed to
apply the Reidemeister moves to  a shadow diagram (without regard to the types of crossings) then the diagram
can be reduced to a disjoint union of circles. This reduction is no longer true in the presence of virtual
crossings.  

\smallbreak

$$\vbox{\picill2inby1.2in(virtual39)  }$$ 

\begin{center}
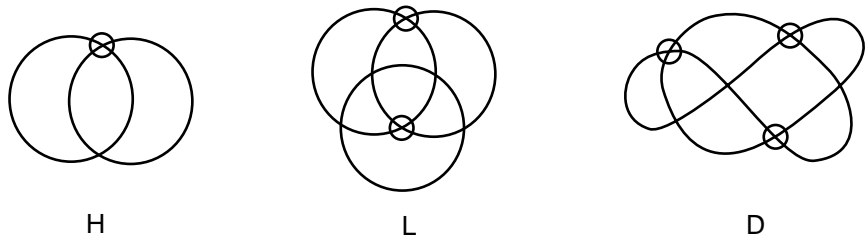
{\bf Figure 39 -- Examples of  Flat Knots and Links} 
\end{center}

More precisely, let a {\it flat virtual diagram} be a diagram with flat  crossings  and virtual  crossings. Two
flat virtual diagrams are {\em equivalent} if there is a sequence of {\it flat virtual Reidemeister moves} taking
one to the other. These are moves as shown in Figure 2, but with flat crossings in place of  classical
crossings. Note that in the category of flat virtuals there is only one forbidden move. 
 Detour moves as in Figure 2C are available only for virtual crossings with respect to flat crossings and not
the other way around.  The study of flat virtual knots and links was initiated in
\cite{VKT}. The category of flat virtual knots is identical in structure to what are called {\em virtual
strings} by V. Turaev  in \cite{TURAEV}. 

\smallbreak
Figure 39 illustrates  flat virtual links $H$ and $L$ and a flat virtual knot $D.$ The link $H$ cannot be
undone in the flat  category because it has an odd number of virtual crossings between its two
components and each flat virtual Reidemeister move preserves the parity of the number of virtual
crossings between components.    The diagram $D$ is shown to be a 
non-trivial flat virtual knot using the filamentation invariant, see \cite{HRK}. The
diagram $L$ is also a non-trivial flat diagram. Note that it comes apart at once if we allow the
forbidden move. 
\smallbreak

Just as virtual knots and links can be interpreted via stabilized embeddings of curves in thickened
surfaces, flat virtuals can be interpreted as stabilized immersions of curves in surfaces (no thickening
required). See \cite{KADOKAMI} for applications of this point of view. Similarly, flat virtual links and braids
have ribbon surface interpretations. In Figure~40 we illustrate the mixed RIII move and its local ribbon surface
embedding. Note the stark difference here between the virtual crossing structure and the immersion structure of 
the flat  crossings.

\smallbreak

$$\vbox{\picill3.5inby.9in(virtual40) }$$ 

\begin{center}
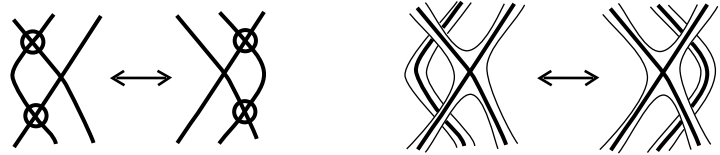
{\bf Figure 40 -- Flat Version of the Detour Move} 
\end{center}

We shall say that a virtual diagram {\em overlies} a flat diagram if the virtual diagram is obtained 
from the flat diagram by choosing a crossing type for each flat crossing in the virtual diagram. To each
virtual diagram $K$ there is an associated flat diagram $F(K)$ that is obtained by forgetting the extra
structure at the classical crossings in $K.$ Note that if $K$ is equivalent to $K'$ as virtual diagrams,
then $F(K)$ is equivalent to $F(K')$ as flat virtual diagrams. Thus, if we can show that $F(K)$ is not
reducible to a disjoint union of circles, then it will follow that $K$ is a non-trivial virtual link. 

\bigbreak

The flat virtual braids were introduced in \cite{SVKT}. As with the virtual braids, the set of flat
virtual braids on $n$ strands forms a group, the {\it flat virtual braid group,}  denoted $FV_n.$ 
The generators of $FV_n$ are the virtual crossings $v_i$ and the flat crossings
$c_i,$ such that $c_i^2 = 1$.  Both, flat crossings and  virtual crossings  represent
geometrically the generators of the symmetric group $S_n.$  But the {\it mixed relation} between them: 

$$v_{i} c_{i+1} v_{i}  =  v_{i+1} c_i v_{i+1}$$

\noindent is not symmetric (see Figure 40).  $FV_n$ is a quotient of the virtual braid group
$VB_{n}$ modulo the relations ${\sigma_i}^2  = 1$ for all $i.$ 
Thus,  $FV_n$ is the free product of two
copies of $S_n,$ modulo the set of mixed relations. Note that $FV_2 = S_2 * S_2$ (no extra
relations), and it is infinite.
\bigbreak

From the above, the flat virtual braids are the appropriate theory of braids for the category of virtual
strings. Every virtual string is the closure of a flat virtual braid. 
 In order to obtain a  Markov theorem for flat virtual braids, we only need to forget, in our study of
virtuals and the definitions of the virtual $L$--moves, the distinction between over and under crossings. The
presence of the flat forbidden move gives rise to the {\it flat threaded  $L_v$--moves}, left and right,  the
analogues of the  threaded $L_v$--moves. Figures~11, 12 and 13 provide illustrations, if we substitute the real crossings
by flat ones. Thus, we have the following results.

\begin{th}[$L$--move Markov Theorem for flat virtuals] Two oriented  flat virtual links are isotopic if and
only if any two corresponding flat virtual braids differ by flat virtual braid isotopy and a finite sequence
of the following moves or their inverses:
\begin{itemize}
\item[(i)]Flat  conjugation 
\item[(ii)]Right virtual $L_v$--moves  
\item[(iii)]Right flat $L_v$--moves  
\item[(iv)]Right and left flat threaded  $L_v$--moves.
\end{itemize} 
\end{th}

\begin{th}[Algebraic Markov Theorem for flat virtuals] Two oriented  flat virtual links are isotopic if and
only if any two corresponding flat virtual braids differ by braid relations in $FV_{\infty}$ and a finite
sequence of the following moves or their inverses:

\begin{itemize}
\item[(i)]Virtual and flat conjugation:  \ \ \ \ \ \ \  \, $ v_i \alpha v_i \sim \alpha \sim 
c_i\alpha c_i $
\item[(ii)]Right virtual and flat stabilization:  \  $\alpha v_n \sim \alpha
\sim \alpha c_n $ 
\item[(iii)]Algebraic right flat threading: \ \ \ \ \ \ \ $\alpha \sim \alpha c_n v_{n-1}
c_n $
\item[(iv)]Algebraic left flat threading: \ \ \ \ \ \ \ \ \ $\alpha \sim  \alpha v_n v_{n-1} c_{n-1}
 v_n c_{n-1} v_{n-1} v_n $
\end{itemize} 

\noindent where $\alpha, v_i, c_i \in FV_n$ and  $ v_n, c_n \in FV_{n+1}.$ (Figure 35 provides   
illustrations, substituting the real crossings by flat ones).
\end{th}

\subsection {Welded Links and Unrestricted Virtuals }

Welded  braids were introduced in \cite{FRR}. They satisfy the same isotopy relations as 
the virtuals, but for welded braids  one of the two forbidden moves of Figure~4 is allowed, the move $F_1$, which 
contains an {\it over arc} and one virtual crossing. One can consider welded knots and links as closures of welded
braids. The move $F_1$ can be regarded as a way of detouring sequences of classical  crossings {\it over} welded
crossings. The explanation for the choice of moves lies in the fact that the  move $F_1$ preserves the  combinatorial
fundamental group. This is not true for the other forbidden move $F_2.$ The   {\it welded braid group}   on $n$ strands, 
$WB_n,$ is a quotient of the  virtual braid group, so it can be presented with the same generators and
relations as  $VB_n,$ but with the extra relations:

$$v_{i} \sigma_{i+1} \sigma_{i}  =  \sigma_{i+1} \sigma_i v_{i+1}   \  \  \ (F_1).$$

In order to obtain a  Markov type theorem for welded braids, we only need to consider in our study of virtuals 
the effect of the move $F_1$. The presence of this move makes redundant the under--threaded  $L_v$--moves
since, by the move $F_1$, the thread can be pulled away, reducing the move to a basic
$vL_v$--move. Thus, threading disappears from the theory of welded braids and we have the following  results  (compare
\cite{Ka}).

\begin{th}[$L$--move Markov Theorem for welded knots] Two oriented  welded links are isotopic if and
only if any two corresponding welded braids differ by welded braid isotopy and a finite sequence of the
following moves or their inverses:
\begin{itemize}
\item[(i)]Real  conjugation 
\item[(ii)]Right virtual $L_v$--moves  
\item[(iii)]Right real $L_v$--moves.  
\end{itemize} 
\end{th}

\begin{th}[Algebraic Markov Theorem for welded knots] Two oriented  welded links are isotopic if and
only if any two corresponding virtual braids differ by braid relations in $WB_{\infty}$ and a finite
sequence of the following moves or their inverses:

\begin{itemize}
\item[(i)]Virtual and real conjugation:  \ \ \ \ \ \ \ \ $ v_i \alpha v_i \sim \alpha \sim 
{\sigma_i}^{-1}\alpha \sigma_i $
\item[(ii)]Right virtual and real stabilization:  \  $\alpha v_n \sim \alpha
\sim \alpha \sigma_n^{\pm 1}$ 
\end{itemize} 

\noindent where $\alpha, v_i, \sigma_i \in WB_n$ and  $ v_n, \sigma_n \in WB_{n+1}$ (recall Figure 35 for
illustrations).
\end{th}

This statement of the Markov Theorem for welded braids is equivalent to that of S.~Kamada \cite{Ka}.

\bigbreak

Finally, another quotient of the virtual braid group (and of the welded braid group) is obtained by adding
both types of forbidden moves. We call this the {\it unrestricted virtual braid group,} denoted $UB_n.$ It is
known that any classical knot can be unknotted in the virtual category if we allow both forbidden moves
\cite{KANENOBU, NELSON}. Nevertheless, linking phenomena still remain.  The unrestricted braid group
itself is non trivial, deserving further study.  For  a presentation of $UB_n$ we just add  to the
presentation of $VB_n$ both types of forbidden moves: 

$$v_{i} \sigma_{i+1} \sigma_{i}  =  \sigma_{i+1} \sigma_i v_{i+1}  \  \ (F_1) \ \ \ \mbox{\it and} \ \
\  \sigma_i \sigma_{i+1} v_{i}  =  v_{i+1}  \sigma_{i} \sigma_{i+1}  \  \ (F_2). $$

Then we have the following:

\begin{th}[Algebraic Markov Theorem for unrestricted virtuals] Two oriented  unrestricted virtual links are
isotopic if and only if any two corresponding unrestricted virtual braids differ by braid relations in
$UB_{\infty}$ and a finite sequence of the following moves or their inverses:

\begin{itemize}
\item[(i)]Virtual and real conjugation:  \ \ \ \ \ \ \ \ $ v_i \alpha v_i \sim \alpha \sim 
{\sigma_i}^{-1}\alpha \sigma_i $
\item[(ii)]Right virtual and real stabilization:  \  $\alpha v_n \sim \alpha
\sim \alpha \sigma_n^{\pm 1}$ 
\end{itemize} 

\noindent where $\alpha, v_i, \sigma_i \in UB_n$ and  $ v_n, \sigma_n \in UB_{n+1}.$ 
\end{th}

Note that the moves of the equivalence relations in Theorems~7 and 8 are apparently the same. The difference in the
theory lies in the different structures of the corresponding braid groups.

\section {On Virtual $R$--matrices and Virtual Hecke Algebras}

In this section we illustrate relations on an R-matrix solution to the Yang-Baxter equation  that would allow
an analog of the Markov trace construction to be made for virtual braids. Such a construction leads to
invariants of virtual knots and links, yielding valuable information about  the virtual category. In Figure 41 we
illustrate the apparatus and relations that are needed to construct a Markov trace on braids from an R-matrix
in the classical case. 
\smallbreak

The illustration uses diagrammatic matrix notation. In this notation a matrix or tensor is  represented by a
box or otherwise delineated polygon in the plane with strands emanating from the box, indicating the indices of
the matrix. When a line from one diagrammatic matrix is tied with  
 a line from another, we see an
internal edge in the graphical structure and this is interpreted as a shared index in the matrix
interpretation. Thus at the matrix level one sums over all possible indices that label an internal edge, and
one  takes the products of all the matrix entries concerned. This is an exact generalization of the formula
for matrix multiplication
$$(MN)_{ab} = \sum_{i} M_{ai}N_{ib}$$ where summation is over all indices $i$ relevant to this matrix product.
\smallbreak 

One can conceptualize diagrammatic
matrices by regarding the diagrams as morphisms in a graphical category, and the intepretation as matrix 
multiplication as a functor to a linear algebraic category. The same remarks apply to the well-known Einstein
summation convention where we write $$M_{ai}N_{ib}$$ and interpret the  repeated index as a summation over
all values for $i.$ Here the algebraic notation $M_{ai}N_{ib}$ is in an abstract tensor category of  indexed
algebra with rules for  handling repeated indices. For example, $M_{ai}N_{ib} = M_{aj}N_{jb}$ so long as $j$
is also repeated and $j$ denotes a letter distinct from $a$ and $b.$ The interpretation as summation takes
the abstract tensor category to a linear algebra category. The diagrams are a generalization of the abstract
tensor category. We use this diagrammatic matrix algebra in our illustrations to show the  translation from the
category of  link diagrams and virtual link diagrams to the matrix algebraic formulas that can capture an
invariant of virtual knots and links via the Markov theorem. For example, see Figure 41.

\bigbreak 

$$ \picill3inby4in(virtual41) $$

\begin{center}
{\bf Figure 41 -- R-Matrix Relations}
\end{center}
\vspace{3mm} 

\noindent  The first
diagram at the upper left denotes a matrix $\eta^{i}_{j}$ where the indices are designated by the strands
emanating from the  black disk that is the body of the diagrammatic version of $\eta.$ Crossings are
represented by the matrices $R = (R^{ij}_{kl})$ and
$\overline{R} = R^{-1}.$ These matrices must satisfy a braiding relation that corresponds to the third 
Reidemeister move. At the matrix level this relation is called the Yang-Baxter Equation. We have not
illustrated this relation. The virtual crossings are shown in Figure 42. They are represented by  a matrix $V$
that must also satisfy the Yang-Baxter equation and the detour relations that generate the virtual braid
group. These matrices then generate a tensor representation of the virtual braid group where a generator
acting on the $i$-th and $(i+1)$-st strands receives an $R,
\overline{R}$ or $V,$ according as it is a classical or virtual crossing, and all the other strands receive an 
identity matrix. Given a virtual braid $\beta,$ let $\rho(\beta)$ denote this  representation applied to
$\beta.$ 
\smallbreak

Now return to Figure 41. Note that we define a trace-function on braids by the formula
$$tr(\beta) = trace(\eta^{\otimes n} \rho(\beta)).$$
Here {\it trace} denotes the usual trace of a matrix.
This trace formula is indicated diagrammatically by the figure in the box to its immediate left.
In order for $tr(\beta)$ to be constructed (after normalization) as a virtual link invariant, we need:  
\begin{enumerate}
\item $tr(\beta \gamma) = tr(\gamma \beta),$ for any braids $\beta$ and $\gamma,$
\item $tr(\beta)$ should either be invariant or it should multiply by a constant under stabilization  moves and
under--threaded moves.
\end{enumerate}

 In Figure 41 we have indicated $tr(\beta)$ to multiply by $\alpha$ under 
right positive classical stabilization and by $\alpha^{-1}$ under right negative classical stabilization.  In
Figure 42 we have indicated right virtual stabilization invariance. In 
Figure 42 we also illustrate the diagrammatics of a right under--threaded move.  Note that these stabilization
equations all involve the matrix $\eta.$ Appropriate choices of the solutions to the Yang-Baxter equation and 
the matrix $\eta$ can, in principle, lead to invariants of both classical and virtual knots and links. One
obtains a normalized invariant
$Invar(b)$ by the formula $$Invar(b) = \alpha^{-w(b)}tr(b)$$ where $w(b)$ is the sum of the signs of the 
exponents of the classical braid generators in  an expression for the braid $b.$
\smallbreak

One case is worth mentioning here explicitly. Suppose that $\eta$ and $R$ yield an invariant of classical 
braids (of which there are many, including the Jones polynomial and specializations of the homflypt
polynomial). Then we can take $V$ (as a linear mapping) to be the permutation
$V(x \otimes y) = y \otimes x.$ Under these conditions $tr(b)$ will satisfy classical
stabilization, but will not necessarily satisfy virtual stabilization. We call such invariants {\it virtual 
rotational invariants}. They are interesting in their own right. It is a subtle matter to obtain full
virtual invariants, but there are examples, including the Jones polynomial itself. 

\bigbreak 

$$ \picill3inby4.8in(virtual42) $$

\begin{center}
{\bf Figure 42 -- Virtual R-Matrix Relations}
\end{center}
\vspace{3mm}

Theorem 3 opens up yet another possibility to construct invariants of virtual links using  algebraic means. 
Namely, to study quotients of the virtual braid group algebra and try to construct on them linear
Markov--type traces. Then, to apply appropriate normalizations yielding virtual link invariants.  Taking the
lead of Jones's construction  \cite{J} of the homflypt (2-variable Jones) polynomial we define ${\cal
VH}_n(q)$, the {\it virtual Hecke algebra} as the quotient of the virtual braid group algebra  $\ZZ [q^{\pm
1}]VB_n$ by factoring out the quadratic relations:  
\[ \sigma_i^2=(q-1)\sigma_i+q. \]  

Let $g_1, \ldots, g_n, v_1, \ldots, v_n$ be the generators of  ${\cal VH}_{n+1}(q).$  A {\it virtual
Markov trace} is defined to be a linear  function $tr$ on $\bigcup_{n=1}^{\infty} {\cal VH}_n(q)$ which 
supports the real and virtual Markov properties. 
 More precisely, we require the trace $tr$ to satisfy the rules: 

\[\begin{array}{ll} 
1) & tr(ab)=tr(ba)     \\  
2) & tr(1)=1          \ \mbox{for all }  {\cal VH}_n(q) \\ 
3) & tr(ag_n)=z\, tr(a)  \\ 
4) & tr(av_n)=s\, tr(a)  \\ 
5) & tr(ag_n^{-1}v_{n-1} g_n^{+1})=r\, tr(a)  \\ 
6)  & tr(a v_n v_{n-1} g_{n-1}^{+1} v_n g_{n-1}^{-1} v_{n-1} v_n )=k\, tr(a)  \\ 
\end{array}\]

\noindent for  $ a,b \in {\cal VH}_n(q)$ and  $z, s, r, k$  independent variables in $\ZZ\, [q^{\pm 1}]$. 
Finally, we  normalize $tr$ appropriately in order to obtain an invariant of virtual links.
\smallbreak

We will pursue these
matters of $R$-matrix invariants of virtual braids and virtual Hecke algebras in a subsequent paper. 
\bigbreak

 \noindent {\bf Acknowledgments.} \ We are happy to mention that the paper of S. Kamada \cite{Ka}
has been for us a source of inspiration. We also thank the referee for remarks that enabled us to sharpen the statements
of the main theorems. It also gives us  pleasure to acknowledge a list of
places and meetings where we worked on these matters. These are, in chronological order: Bedlewo, Frankfurt airport,
Oberwolfach, Caen, Athens, Vancouver, Albuquerque, Chicago, Oberwolfach, Chicago, Athens. The last stay at Oberwolfach,
where the paper was completed, was due to the Research In Pairs program, which we
acknowledge gratefully. It gives the first author great pleasure to acknowledge support from  NSF Grant
DMS-0245588.

\bigbreak

\noindent {\sc L.H. Kauffman: Department of Mathematics, Statistics and
Computer Science, University of Illinois at Chicago, 851 South Morgan St., Chicago IL 60607-7045, USA. }

 \vspace{.1in}
 
\noindent {\sc S. Lambropoulou:  Department of Mathematics,  National Technical University of Athens,
Zografou campus, GR-157 80 Athens, Greece.}

\vspace{.1in}
\noindent {\sc E-mails:} \ {\tt 
kauffman@uic.edu  \ \ \ \ \ \ \ \ \ \ \ \ \ \ \ \ sofia@math.ntua.gr}

\noindent {\sc URLs:} \ {\tt 
www.math.uic.edu/$\tilde{~}$kauffman  \ \ \ \ \ \ \ \ http://www.math.ntua.gr/$\tilde{~}$sofia/}


\end{document}